\documentclass[a4paper,12pt]{article}  
\RequirePackage[OT1]{fontenc}  
\usepackage[natbib,noinfoline]{imsart}  
\usepackage[OT1]{fontenc}  
\usepackage[english]{babel}  
\usepackage[utf8]{inputenc}   
\usepackage{latexsym,amssymb,amsmath,epsfig,psfrag}  
  
\def\PP{\mathbb{P}} 
\def\RR{\mathbb{R}} 
\def\EE{\mathbb{E}} 
\def\NN{\mathbb{N}} 
\def\ZZ{\mathbb{Z}} 
\def\SS{\mathbb{S}}
\def\Ent{\mathsf{Ent}} 
\def\LL{\mathbf{L}}
\def\Var{\mathsf{Var}}
\def\II{\mbox{ 1\hskip -.29em I}}
\newcommand{\sfrac}[2]{\kern.1em
        \raise.5ex\hbox{$#1$}\kern-.1em
        /\kern-.15em\lower.25ex\hbox{$#2$}}
\def\bdes{\begin{description}}
\def\edes{\end{description}}

\def\iti{\item[(i)]}
\def\itii{\item[(ii)]}
\def\itiii{\item[(iii)]}
\def\itiv{\item[(iv)]}
\def\itv{\item[(v)]}

\newtheorem{defi}{Definition}[section]
\newtheorem{lemm}[defi]{Lemma}
\newtheorem{prop}[defi]{Proposition}
\newtheorem{coro}[defi]{Corollary}
\newtheorem{theo}[defi]{Theorem}
\newtheorem{exem}{Example}
\newtheorem{rem}{Remark}
\newenvironment{dem}{\vskip 2mm\noindent {\it Proof} :}
                    {\hfill $\square$ \vskip 2mm \noindent}

\def\bex{\begin{exem} \em }
\def\eex{\end{exem} }
\def\brem{\begin{rem} \em}
\def\erem{\end{rem} }

\def\mylabel#1{\label{#1}}

\begin{document}
\begin{frontmatter}

\title{Exponential concentration for First Passage Percolation
  through modified Poincar\'e inequalities} 

\runtitle{Modified Poincaré inequalities \& FPP}
\thanks{We acknowledge financial support from the Swiss National Science
  Foundation grants 200021-1036251/1 and 200020-112316/1.} 
\author{\fnms{Michel} \snm{Benaïm} \ead[label=e2]{michel.benaim@unine.ch}}
\address{Institut de Mathématiques\\ Universit\'e de Neuchâtel,\\ 11 rue Emile
  Argand,\\ 2000 Neuchâtel, SUISSE\\
\printead{e2}}
\author{\fnms{Rapha\"el} \snm{Rossignol} \ead[label=e1]{raphael.rossignol@unine.ch}}
\address{Institut de Mathématiques\\ Universit\'e de Neuchâtel,\\ 11 rue Emile
  Argand,\\ 2000 Neuchâtel, SUISSE\\
\printead{e1}}


\runauthor{M. Benaïm, R. Rossignol}

\begin{abstract}:~ We provide a new exponential concentration inequality for
  First Passage Percolation valid for a wide class of edge times
  distributions. This improves and extends a result by Benjamini, Kalai and Schramm
  \cite{BenjaminiKalaiSchramm03} which gave a variance bound for Bernoulli
  edge times. Our approach is based on some functional
  inequalities extending the work of Rossignol \cite{Rossignol06}, Falik and
  Samorodnitsky \cite{FalikSamorodnitsky06}.\vspace*{3mm}

\noindent {\bf Résumé}:~ On obtient une nouvelle inégalité de concentration
  exponentielle pour la percolation de premier passage, valable pour une large
  classe de distributions des temps d'arêtes. Ceci améliore et étend un
  résultat de  Benjamini, Kalai et Schramm
  \cite{BenjaminiKalaiSchramm03} qui donnait une borne sur la variance pour
  des temps d'arêtes suivant une loi de Bernoulli. Notre approche se fonde sur des inégalités
  fonctionnelles étendant les travaux de Rossignol \cite{Rossignol06}, Falik et
  Samorodnitsky \cite{FalikSamorodnitsky06}.
\end{abstract}



\end{frontmatter}

\section{Introduction}
\mylabel{sec:intro}
First Passage Percolation was introduced by Hammersley and Welsh
\cite{HammersleyWelsh65} to model the
flow of a fluid in a randomly porous material (see \cite{Howard04} for a recent
account on the subject). We will consider
the following model of First Passage Percolation in $\ZZ^d$, where $d\geq 2$ is an
integer. Let $E=E\left(\ZZ^d\right)$ denote the set of edges in $\ZZ^d$. The
passage time of the fluid through the edge $e$ is denoted by
$x_e$ and is supposed to be nonnegative. Randomness of the porosity is given
by a product probability measure on $\RR_+^E$. Thus,
$\RR_+^E$ is equipped with the  measure $\mu=\nu^{\otimes E}$,
where $\nu$ is a probability measure on $\RR_+$ according to which each
passage time is distributed, independently from the others. If $u,v$ are two vertices of
$\ZZ^d$, the notation $\alpha:\{u,v\}$ means that $\alpha$ is a path with
end points $u$ and $v$. When $x\in\RR_+^E$, $d_x(u,v)$ denotes the first passage
time, or equivalently the distance from $u$ to $v$ in the metric induced by $x$,
$$d_x(u,v)=\inf_{\alpha:\{u,v\}}\sum_{e\in\alpha}x_e\;.$$
The study of $d_x(0,nu)$ when $n$ is an integer
which goes to infinity is of central importance. Kingman's subadditive
ergodic theorem implies the existence, for each fixed $u$, of a ``time
constant'' $t(u)$ such that:
$$\frac{d_x(0,nu)}{n}\xrightarrow[n\rightarrow +\infty]{\nu-a.s} t(u)\;.$$
It is known (see Kesten \cite{Kesten84}, p.127 and 129) that if $\nu(\{0\})$ is
strictly smaller than the critical probability for Bernoulli bond percolation
on $\ZZ^d$, then $t(u)$ is  positive for every $u$ distinct from the
origin. Under such an  assumption, one can say that the random variable $d_x(0,nu)$
is located around $nt(u)$, which is of order $O(|nu|)$, where we denote by $|.|$ the $L^1$-norm of vertices in $\ZZ^d$. In this paper, we
are interested in the fluctuations of this quantity. Precisely, we define, for
any vertex $v$,
$$\forall x\in\RR_+^E,\; f_v(x)=d_x(0,v)\;.$$
It is widely believed that the fluctuations of $f_v$ are of order
$|v|^{1/3}$ when $d=2$. Apart from some predictions made by physicists, this faith relies on recent results for related growth models \cite{BaikDeiftJohansson99,Johansson00,Johansson00phys}. Until recently, the best
results rigourously obtained for the fluctuations of $f_v$ were some
moderate deviation estimates of order $O(|v|^{1/2})$ (see
\cite{Kesten93,Talagrand95}). In 1993, Kesten \cite{Kesten93} proved that,
$$\Var_\mu(f_v)=O(|v|)\;,$$
provided
$\nu$ admits a finite second order moment. If, furthermore, $\nu$ admits a finite moment of exponential order, there exist two
constants $C_1$ and $C_2$ such that for any $t\leq |v|$,
\begin{equation}
\label{ineqKestenexpo}
\nu(|f_v-\EE(f_v)|>t\sqrt{|v|})\leq C_1e^{-C_2t}\;.
\end{equation}
Later, Talagrand improved the right-hand side of the above inequality to
$\exp (-C_2t^2)$. In 2003, Benjamini, Kalai and Schramm
\cite{BenjaminiKalaiSchramm03} proved that for Bernoulli
edge times bounded away from 0, the variance of $f_v$ is of order $O(|v|/\log|v|)$, and therefore,
the fluctuations are of order  $O(|v|^{1/2}/(\log|v|)^{1/2})$. 

It is natural to ask whether the work of Benjamini, Kalai and Schramm
\cite{BenjaminiKalaiSchramm03} can be extended to other distributions, notably
continuous distributions which are not bounded away from zero. This has been done
in a preliminary version of the present paper \cite{BenaimRossignolarxiv06} by
extending the tools of \cite{BenjaminiKalaiSchramm03}, namely a modified
Poincaré inequality due to Talagrand. It is also natural, and even more desirable, to try to improve
the result of Benjamini, Kalai and Schramm \cite{BenjaminiKalaiSchramm03} into
an exponential inequality in the spirit of (\ref{ineqKestenexpo}), with
$\sqrt{|v|/\log |v|}$ instead of $\sqrt{|v|}$. In this article we show that some different modified Poincaré
inequalities arising from the context of ``threshold phenomena'' for Boolean
functions (see \cite{Rossignol06,FalikSamorodnitsky06}) may be used
successfully instead of Talagrand-type inequalities from
\cite{Talagrand94a,BenaimRossignolarxiv06}. This is the main result of this
paper, stated in Theorem \ref{theoFPPexpo}.  Whereas we focused on the percolation setting, the argument is fairly
general and we present also an abstract exponential concentration result,
Theorem \ref{thmconcentrationgenerale}, which is very likely to have
applications outside the setting of percolation.


This article is organized as follows. In Section \ref{sec:logsobfkbis}, we
extend the modified Poincaré inequalities of Falik and Samorodnitsky
\cite{FalikSamorodnitsky06} to non-Bernoulli and countable settings where
logarithmic Sobolev inequalities are available, notably
to a countable product of Gaussian measures. Section \ref{sec:extension} is devoted to the obtention
of similar inequalities for other continuous measures by a simple mean of
change of variable. In section \ref{sec:concentration}, we show how to deduce
new general exponential concentration bounds from the modified Poincaré inequalities of section
\ref{sec:logsobfkbis}. This allows us to obtain in section \ref{sec:FPP} an exponential version of the
bound of Benjamini et
al. in some continuous and discrete settings.

\paragraph{Notation}
Given a probability space $(\mathbb{X}, {\cal X}, \mu)$ and a real valued
measurable function $f$ defined on $\mathbb{X}$ we let  $$\|f\|_{p,\mu} =
\left(\int |f|^p\;d\mu\right)^{1/p} \in [0, \infty]\;,$$
and $L^p(\mu)$ denote the set of $f$ such that $\|f\|_{p,\mu} < \infty.$ 
The {\em mean} of $f \in L^1(\mu)$ is denoted
$$\EE_{\mu}(f)  = \int f d\mu\;,$$
the {\em variance} of $f \in L^2(\mu)$ is 
$$\Var_{\mu}(f) = \|f - \EE_{\mu}(f)\|_{2,\mu}^2\;,$$
and the {\em entropy} of any positive measurable function $f$ is
$$\Ent_{\mu}(f) = \int f\log f\;d\mu-\int f\;d\mu\log\int f\;d\mu\;.$$
When the choice of $\mu$ is unambiguous we may write $\|f\|_p$ (respectively
$L^p$, $\EE(f)$, $\Var(f)$ and $\Ent(f)$)    for $\|f\|_{p,\mu}$
(respectively  $L^p(\mu)$, $\EE_{\mu}(f)$, $\Var_{\mu}(f)$ and $\Ent_{\mu}(f)$). 

\section{Logarithmic Sobolev and modified Poincar\'e inequalities on $\RR^{\NN}$}
\mylabel{sec:logsobfkbis}

The relevance of a ``modified Poincaré inequality'' due to Talagrand
\cite{Talagrand94a} in the context of First Passage Percolation was shown by
Benjamini, Kalai and Schramm \cite{BenjaminiKalaiSchramm03}. Let us explain
this point a little bit more. A classical Poincaré inequality has the following form:
$$\Var_\mu(f)\leq C\mathcal{E}_\mu(f)\;,$$
where $C$ is a constant, and $\mathcal{E}_\mu(f)$ is an ``energy'' of $f$, that
is, usually, the mean against $\mu$ of the square of some kind of
gradient. There is a good theory for this in the context of Markov semi-groups
(see \cite{thesardsledoux,Bakry02tata} for instance). By ``modified Poincaré
inequality'', we mean a functional inequality which improves upon the
classical Poincaré inequality for a certain class of functions $f$. This is
usually achieved through a hypercontrativity property (see
\cite{BenaimRossignolarxiv06} and Ledoux \cite{LedouxJerusalem}).

In this section, we will show how to build a modified Poincaré inequality on
a product of probability spaces each of which satisfies a Sobolev logarithmic
inequality. This approach was initiated independently by Rossignol
\cite{Rossignol06}, Falik and Samorodnitsky \cite{FalikSamorodnitsky06} in the
Bernoulli setting. 

We shall need some notation for tensorisation. Suppose that we are
given a countable collection of probability spaces
$\left(\mathbb{X}_i,\mathcal{X}_i,\mu_i\right)_{i\in I}$.   If $i$ belongs to  $I$, and $x^{-i}$ is an element of
$\Pi_{\substack{j\in I\\ j\not=i}}\mathbb{X}_j$, then for every $x_i$ in $\mathbb{X}_i$, we denote by
$(x^{-i},x_i)$ the element $x$ of $\Pi_{j\in I}\mathbb{X}_j$. 
For every function $f$ from $\Pi_{i\in I}\mathbb{X}_i$ to $\RR$,
every $j\in I$, and every $x^{-j}$ in $\Pi_{i\not =j}\mathbb{X}_i$,  we denote by $f_{x^{-j}}$ the function from
$\mathbb{X}_j$ to $\RR$ obtained from $f$ by keeping $x^{-j}$ fixed:
$$\forall x_j\in\mathbb{X}_j,\;
f_{x^{-j}}(x_j)=f(x)\;.$$
Now, suppose that we are given a collection
$(\mathcal{A}_i)_{i\in I}$ of  linear subspaces, $\mathcal{A}_i\subset L^2(\mathbb{X}_i,\mu_i)$, containing the constant
functions. Then, we introduce
$$\mathcal{A}^I=\{f\in L^2(\Pi_{i\in I}\mathbb{X}_i,\;\otimes_{i\in
  I}\mu_i)\mbox{ s.t. }\forall j\in I,\;f_{x^{-j}}\in
\mathcal{A}_j\mbox{ for }\otimes_{i\not =j}\mu_i\mbox{-a.e }x^{-j}\}\;.$$
An operator $R_j$ from $\mathcal{A}_j$ to $L^2(\mathbb{X}_j,\mu_j)$ is naturally extended on $\mathcal{A}^I$, `` acting only on coordinate $j$ '':
$$\forall f\in \mathcal{A}^I,\;\forall x\in\Pi_{i\in
  I}\mathbb{X}_i,\;R_j(f)(x):=R_j(f_{x^{-j}})(x_j)\;.$$

\begin{prop}
\label{proppoincaremod}
 Let $\left(\mathbb{X}_i,\mathcal{X}_i,\mu_i\right)_{i\in I}$, be a sequence
 of probability
 spaces. Let $(\mathcal{A}_i)_{i\in I}$ be a collection of sets such that for every $i$, $\mathcal{A}_i$ is
a linear subspace of $L^2(\mathbb{X}_i,\mu_i)$ which contains the constant
functions. Suppose that for every $i$ in $I$, $\mu_i$ satisfies a logarithmic Sobolev inequality of the
 following form:
$$\forall f\in \mathcal{A}_i,\; \Ent_{\mu_i}(f^2)\leq \EE_{\mu_i}\left(R_i(f)^2\right)\;,$$
where $R_i$ is a linear operator from  $\mathcal{A}_i$ to $L^2(\mu_i)$ with value zero on any
constant function. Furthermore, suppose that the following commutation property holds:
$$\forall i,\;\forall f\in \mathcal{A}^{I},\;\int
f\;d\otimes_{j\not=i}\mu_j\in\mathcal{A}_i\mbox{ and }R_i\left(\int f\;d\otimes_{j\not=i}\mu_j\right)=\int
R_i(f)\;d\otimes_{j\not=i}\mu_j\;.$$
Then, $\mu^I=\otimes_{i\in I}\mu_i$ satisfies the
following modified Poincaré inequality:
$$\forall f\in
\mathcal{A}^{\NN},\;\Var_{\mu^I}(f)\log\frac{\Var_{\mu^I}(f)}{\sum_{i\in
    I}\left\|\Delta_i f\right\|_{\mu^I,1}^2}\leq
\sum_{i\in I}\EE_{\mu^I}\left(R_i(f)^2\right)\;,$$
where $\Delta_i$ is the following operator on $L^2(\mu^I)$:
$$\forall f\in L^2(\mu^I),\;\Delta_if=f-\int f\;d\mu_i\;.$$
\end{prop}

\begin{dem}
To shorten
  the notations, we shall write $\mu$ instead of $\mu^I$.

First, suppose that $I$ is finite, $I=\{1,\ldots ,n\}$. The tensorisation
property of the entropy (see \cite{thesardsledoux} or \cite{Ledoux96}, Proposition 5.6 p.98 for instance) states that for every positive measurable function $g$,
$$\Ent_{\mu}(g)\leq
\sum_{i=1}^n\EE_\mu(\Ent_{\mu_i}(g))\;.$$
Thus, the logarithmic Sobolev inequalities for each $\mu_i$ imply that:
\begin{equation}
\label{eqlogsobproduit}
\forall g\in \mathcal{A}^n,\; \Ent_{\mu}(g^2)\leq
\sum_{i=1}^n\EE_\mu\left(R_i(g)^2\right)\;.
\end{equation}
Now, let $f$ be a function in $\mathcal{A}^n$. Following Rossignol \cite{Rossignol06}, Falik and Samorodnitsky \cite{FalikSamorodnitsky06}, we write $f-\EE_\mu(f)$ as a sum of
martingale increments, and apply the logarithmic Sobolev inequality
(\ref{eqlogsobproduit}) to each increment:
\begin{equation}
\label{eqlogsobincrements}
\sum_{j=1}^n\Ent_{\mu}(V_j^2)\leq
\sum_{j=1}^n\sum_{i=1}^n\EE_\mu\left(R_i(V_j)^2\right)\;,
\end{equation}
where
$$f-\EE_\mu(f)=\sum_{j=1}^nV_j\;,$$
and
$$V_j=\int f \;d\mu_1\otimes\ldots \otimes d\mu_{j-1} - \int f
\;d\mu_1\otimes\ldots \otimes d\mu_{j}\;=\int \Delta_jf\;d\mu_1\otimes\ldots
\otimes d\mu_{j-1}.$$
The following inequality, which is a clever application of Jensen's inequality, is shown in Falik and Samorodnitsky
\cite{FalikSamorodnitsky06} and is cleaner than the corresponding one in
Rossignol \cite{Rossignol06}:
$$\sum_{j=1}^n\Ent_{\mu}(V_j^2)\geq
\Var_\mu(f)\log\frac{\Var_\mu(f)}{\sum_{j=1}^n\left\|V_j\right\|_{\mu,1}^2}\;.$$
Jensen's inequality implies that:
\begin{equation}
\label{eqminoentropie}
\sum_{j=1}^n\Ent_{\mu}(V_j^2)\geq
\Var_\mu(f)\log\frac{\Var_\mu(f)}{\sum_{j=1}^n\left\|\Delta_j
    f\right\|_{\mu,1}^2}\;.
\end{equation}
On the other hand, for every $i$, the term $\EE_\mu\left(R_i(g)^2\right)$ in
  (\ref{eqlogsobproduit}) is called an ``energy'' for $g$, and we claim that the sum of the energies of the increments of
$f$ equals the energy of $f$:
\begin{equation}
\label{eqsommeenergies}
\forall i\in\{1,\ldots ,n\},\;\sum_{j=1}^n\EE_\mu\left(R_i(V_j)^2\right)=\EE_\mu\left(R_i(f)^2\right)\;.
\end{equation}
Indeed, since $R_i$ is linear, using the commutation hypothesis, and the fact
that $R_i(f)$
is zero on any function $f$ which is constant on coordinate $i$, we get:
\begin{eqnarray*}
&&\forall i< j,\; R_i(V_j)=0\;,\\
&&R_i(V_i)=\int R_i(f)\;d\mu_1\otimes\ldots \otimes d\mu_{i-1}\;,
\end{eqnarray*}
and
$$\forall i >j,\;R_i(V_j)=\int R_i(f)\;d\mu_1\otimes\ldots \otimes d\mu_{j-1}-
\int R_i(f)\;d\mu_1\otimes\ldots \otimes d\mu_{j}\;.$$
Therefore,
\begin{eqnarray*}
\sum_{j=1}^n\EE_\mu\left(R_i(V_j)^2\right)&=&\sum_{j=1}^{i-1}\EE_\mu\left(R_i(V_j)^2\right)+\EE_\mu\left(R_i(V_i)^2\right)\;,\\
&=&\EE_\mu\left(R_i(f)^2\right)\;.
\end{eqnarray*}
Now, claim (\ref{eqsommeenergies}) is proved and the result follows from
(\ref{eqsommeenergies}), (\ref{eqminoentropie}) and
(\ref{eqlogsobincrements}), at least when $I$ is finite.

Now, suppose that $I$ is strictly countable, let us say $I=\NN$, and let $\mathcal{F}_n$ be the $\sigma$-algebra generated by the first $n$
coordinate functions in $\RR^\NN$. Let $f \in \mathcal{A}^{\NN}$ and $f_n = \EE\left(f|\mathcal{F}_n\right)$ be the
conditional expectation of $f$ with respect to $\mathcal{F}_n$. Then, the
commutation property tells us that $f_n$ belongs to $\mathcal{A}^{\NN}$, and
$R_i(f_n)=\EE\left(R_i(f)|\mathcal{F}_n\right)$. Therefore, we can apply the
first part of Proposition \ref{proppoincaremod}, the one that we just proved:
$$\Var_{\mu^n}(f_n)\log\frac{\Var_{\mu^n}(f_n)}{\sum_{i=1}^n\left\|\Delta_i f_n\right\|_{\mu^n,1}^2}\leq
\sum_{i=1}^n\EE_{\mu^n}\left(R_i(f_n)^2\right)\;.$$
This may be written as:
$$\Var_{\mu^{\NN}}(f_n)\log\frac{\Var_{\mu^{\NN}}(f_n)}{\sum_{i=1}^n\left\|\Delta_i f_n\right\|_{\mu^n,1}^2}\leq
\sum_{i=1}^n\EE_{\mu^\NN}\left(\EE\left(R_i(f)|\mathcal{F}_n\right)^2\right)\;.$$
Obviously, $\Delta_i f_n=\EE\left(\Delta_i(f)|\mathcal{F}_n\right)$. Therefore,
Jensen's inequality implies:
$$\Var_{\mu^{\NN}}(f_n)\log\frac{\Var_{\mu^{\NN}}(f_n)}{\sum_{i=1}^n\left\|\Delta_i f\right\|_{\mu^n,1}^2}\leq
\sum_{i=1}^n\EE_{\mu^\NN}\left((R_i(f))^2\right)\;.$$
Of course, $f_n$ converges to $f$ in $L^2(\mu^\NN)$, and we may let $n$ tend
to infinity in the last inequality to get the desired result.
\end{dem}

\begin{rem}
Actually, a logarithmic Sobolev inequality associated to a probability measure $\mu_i$
which is  reversible with respect to an operator $\LL$ may always be written in the form of
Proposition \ref{proppoincaremod}. Indeed, such an inequality may be written as:
$$\forall f\in \mathcal{A}_i,\; \Ent_{\mu_i}(f^2)\leq c
\EE_{\mu_i}\left(-f\LL f\right)\;,$$
where $c$ is a positive constant. Since $\LL$ is a self-adjoint operator in $L^2(\mu_i)$, it admits a spectral
representation (see Yosida \cite{Yosida} p.313). It is easy to show that its eigenvalues are non-negative (see,
for instance Bakry \cite{Bakry02tata} p.7). The spectral decomposition of $-\LL$ may
therefore be written as:
$$-\LL=\int_0^{\infty}\lambda\; dE(\lambda)\;,$$
and a suitable candidate for $R_i$ may be deduced from it:
\begin{equation}
\label{eqracine}
R_i=\int_0^{\infty}\sqrt{c\lambda}\; dE(\lambda)\;.
\end{equation}
Nevertheless, in the applications which follow, it is essential that the
operator $R_i$ is nice enough to allow the quantity 
$\sum_{i=1}^n R_i(f)^2$ to be easily controlled, and the
one given in (\ref{eqracine}) may not be appropriate for this. Another
candidate, which we shall see to be the right one for certain continuous
probability measures, is the square root of the ``carré du champ''
operator $\Gamma^{1/2}(f,f)$, where:
$$\Gamma(f,g)=\frac{1}{2}\left(\LL fg-f\LL g-g\LL f\right)\;.$$
But in the discrete case, this is not the most natural choice. Therefore, we prefer
not to try to generalize any longer, and rather give some examples.
\end{rem}

\subsection{Examples}
\label{subsec:examples}

Not surprisingly, we start to illustrate Proposition \ref{proppoincaremod} with the Bernoulli and Gaussian cases. Our choice
to present them ``mixed'' might look a little weird at first sight, but this will prove
to be useful in the percolation context (see section \ref{sec:FPP}). 

\begin{exem} {\bf The Bernoulli and Gaussian cases.}

{\rm We let $$\beta_p = (1-p)\delta_0 + p\delta_1$$ be the
Bernoulli measure with parameter $p$ on $\{0,1\}$. If $p$ belongs to $]0,1[$, $\beta_p$ (see for instance Saloff-Coste
\cite{Saloff} Theorem 2.2.8 p.336, or \cite{thesardsledoux}) satisfies the following logarithmic Sobolev
 inequality: for any function $f$ from $\{0,1\}$ to $\RR$,
$$\Ent_{\beta_p}(f^2)\leq c_{LS}(p)\EE_{\beta_p}\left((\Delta f)^2\right)\;,$$
where 
$$ c_{LS}(p)=\frac{\log p-\log (1-p)}{p-(1-p)}\;,$$
and
$$\Delta f=f-\int f\;d\beta_p\;.$$
If $S$ is a 
countable set, for any $s$ in $S$, let $\mathbb{X}_s$ be a copy of $\{0,1\}$,
$\mathcal{A}_s$ be the set of functions from  $\mathbb{X}_s$ to $\RR$,
 and $R_s$ be the operator $\Delta$ acting on $\mathcal{A}_s$. We denote also
 a product measure $\lambda_p^S$ on $\{0,1\}^S$: $\lambda_p^S =
 \beta_p^{\otimes^S}$.

Now we introduce the Gaussian setting. Let 
$$\gamma(dy) = \frac{1}{\sqrt{2\pi}} e^{-y^2/2} dy$$ denote the  standard
Gaussian measure on $\RR$. A map $f$ is said to be  {\em weakly differentiable} provided there exists a locally integrable function denoted  $f'(x)$ such that
$$\int f'(x)g(x) dx = - \int g'(x) f(x) dy$$
for every smooth function $g : \RR \mapsto \RR$ with compact support. The {\em weighted Sobolev} space $H_1^2\left(\gamma\right)$ is defined to be the space of weakly differentiable functions $f$ on $\RR$ 
such that
$$\|f\|^2_{H_1^2} = \|f\|_2^2 + \left\|f'\right\|_2^2 < \infty.$$
It is well known that $\gamma$ (see for instance Ledoux
\cite{Ledoux01} Theorem 5.1 p.92) satisfies the following logarithmic Sobolev inequality: for any function $f$ in $H_1^2\left(\gamma\right)$,
$$\Ent_{\gamma}(f^2)\leq 2\EE_{\gamma}\left((f'(x))^2\right)\;.$$
For any $i$ in $\NN$, let $\mathbb{X}_i$ be a copy of $\RR$, $\mathcal{A}_i$ a
copy of $H_1^2\left(\gamma\right)$ and $R_i$ the derivation operator on $\mathcal{A}_i$. We let $\gamma^\NN =
\gamma^{\otimes^\NN}$  denote the standard Gaussian measure on $\RR^\NN$. 

The
set $\mathcal{A}^{S\cup \NN}$ thus defined is the so called {\em weighted Sobolev} space
$H_1^2\left(\lambda_p^S\otimes \gamma^\NN\right)$, which contains the functions $f \in L^2\left(\lambda_p^S\otimes \gamma^\NN\right)$ verifying the following condition. For all $i\in \NN,$
there exists a function $h_i$ in $L^2\left(\lambda_p^S\otimes \gamma^\NN\right)$
such that $$-\int_{\RR} g'(y_i)f(x,y)\;dy_i=\int_{\RR} g(y_i)h_i(x,y)\;dy_i\;, \lambda_p^S\otimes \gamma^\NN \; a.s $$
for every smooth function $g :\RR \mapsto \RR$ having compact support.
The function $h_i$ is called the partial derivative of $f$ with respect to $y_i$,
and is denoted by $\frac{\partial f}{\partial
    y_i}.$

Thus, we deduce from Proposition \ref{proppoincaremod} the following result.
\begin{coro}
\label{coroPoincaremodgauss}
For any $p\in]0,1[$, and any $f\in H_1^2\left(\lambda_p^S\otimes \gamma^\NN\right)$,
$$\Var(f)\log\frac{\Var(f)}{\sum_{s\in S}\left\|\Delta_s
    f\right\|_{1}^2+\sum_{i\in\NN}\left\|\Delta_i f\right\|_{1}^2}\leq c_{LS}(p)
\sum_{s\in S}\EE\left((\Delta_s f)^2\right)+2\sum_{i\in\NN}\EE\left(\left(\frac{\partial
    f}{\partial x_i}\right)^2\right)\;.$$
\end{coro}
}
\end{exem}

\begin{exem} {\bf The gamma case (associated to the Laguerre generator).}

{\rm 
We let  
$$\nu_{a,b}(dy) = \frac{b^a}{\Gamma(a)}y^{a-1}e^{-by} \II_{y>0}\;dy$$ denote the
gamma probability measure
with parameters $a$ and $b$. This measure is the invariant distribution of the
Laguerre semi-group, with generator:
$$\LL_{a,b}f(x)=bxf''(bx)-(a-bx)f'(bx)\;.$$
When $a \geq 1/2$ and $b=1$, it can
be easily seen that this generator satisfies the $CD(\rho,\infty)$
curvature inequality:
$$\Gamma_2(f)\geq \frac{1}{2}\Gamma(f)\;.$$
This implies that $\nu_{a,b}$ satisfies the following
logarithmic Sobolev inequality (see Definition 3.1 p.28 and
Theorem 3.2 p.29 in Bakry \cite{Bakry02tata}, see also \cite{thesardsledoux}). For any weakly differentiable
function $f\in L^2(\nu_{a,b})$, if $a\geq 1/2$,
$$\Ent_{\nu_{a,b}}(f^2)\leq
\frac{4}{b}\EE_{\nu_{a,b}}\left((\sqrt{x}f'(x))^2\right)\;.$$
Therefore, we deduce the following result from Proposition \ref{proppoincaremod}.
\begin{coro}
\label{coroPoincaremodgamma}
Suppose that $a\geq 1/2$, $b> 0$, and let ${\RR^+_*}^{\NN}$ be equipped with the product measure
$\nu_{a,b}^{\NN}$. For any weakly differentiable
function $f$ in $L^2(\nu_{a,b}^{\NN})$, define 
$${\nabla}_if(x) =  \frac{\partial f}{\partial x_i}(x) \sqrt{x_i}\;.$$
Suppose that,
$$\forall i\in \NN,\;\nabla_if\in L^2\left(\nu_{a,b}^\NN\right)\;.$$
Then,
$$\Var(f)\log\frac{\Var(f)}{\sum_{i\in\NN}\left\|\Delta_i f\right\|_{1}^2}\leq \frac{4}{b}
\sum_{i\in\NN}\EE\left(\left(\nabla_if\right)^2\right)\;.$$
\end{coro}
Remark that when $a\in ]0,1/2[$, $\nu_{a,b}$ still satisfies a logarithmic Sobolev inequality with a positive constant $C_{a,b}$ instead of $4/b$, but the precise value of $C_{a,b}$ is not known; see \cite{Miclo02}. This gives the analogue of Corollary \ref{coroPoincaremodgamma} for $a\in ]0,1/2[$, with $C_{a,b}$ instead of $4/b$.
}
\end{exem}

\begin{exem} {\bf The uniform case.}

{\rm
We let  
$$\lambda(dy) = \II_{0\leq y\leq1}dy$$ denote the uniform probability measure on $[0,1]$. It is known that $\lambda$ satisfies the following
logarithmic Sobolev inequality (it is a direct consequence of the logarithmic
Sobolev inequality on the circle \cite{EmeryYukich87}). For any weakly differentiable
function $f$ in $L^2(\lambda)$,
$$\Ent_{\lambda}(f^2)\leq \frac{2}{\pi^2}
\EE_{\lambda}\left((f'(x))^2\right)\;.$$
\begin{coro}
\label{coroPoincaremoduniform01}
Let $[0,1]^{\NN}$ be equipped with the product measure $\lambda^{\NN}$. Suppose that,
$$\forall i\in \NN,\;\frac{\partial f}{\partial x_i}\in L^2\left(\lambda^\NN\right)\;.$$
For any weakly differentiable
function $f$ in $L^2(\lambda^{\NN})$,
$$\Var(f)\log\frac{\Var(f)}{\sum_{i\in\NN}\left\|\Delta_i f\right\|_{1}^2}\leq \frac{2}{\pi^2}
\sum_{i\in\NN}\EE\left(\left(\frac{\partial
    f}{\partial x_i}\right)^2\right)\;.$$
\end{coro}
We shall see in section \ref{sec:extension} that $\lambda$ satisfies another logarithmic
Sobolev inequality with an energy whose form ``looks like'' the energy
appearing in the gamma case.
}
\end{exem}

\section{Extension from the Gaussian case to other measures}
\mylabel{sec:extension}
As usual, we can deduce from Corollary \ref{coroPoincaremodgauss} other
inequalities by mean of change of variables. To make this precise, let $\Omega$ be a measurable space and $\Psi : \RR^{\NN} \mapsto \Omega$ a measurable isomorphism  (meaning that $\Psi$ is one to one with $\Psi$ and $\Psi^{-1}$ measurables). Let $\Psi^*\gamma^{\NN}$ denote the image of $\gamma^{\NN}$ by $\Psi.$ That is $\Psi^*\gamma^{\NN}(A) = \gamma^{\NN}(\Psi^{-1}(A)).$ For $g : S \times \Omega \mapsto \RR$ such that $g \circ (Id,\Psi) \in H_1^2(\lambda \otimes \gamma^{\NN}),$ one obviously has
 $$\Var_{\lambda \otimes \Psi^*\gamma^{\NN}}(g) = \Var_{\lambda \otimes \gamma^{\NN}}(g \circ \Psi)$$
and
$$\|\partial_{i,\Psi} g\|_{p,\Psi^*\gamma^{\NN}} = \left\|\frac{\partial g \circ \Psi}{\partial y_i}\right\|_{p,\gamma^{\NN}}$$
where $\partial_{i,\Psi} g$ is defined as
$$\partial_{i,\Psi} g(x,\omega) = \frac{\partial (x \circ \Psi)}{\partial y_i}(q,\Psi^{-1}(\omega))\;.$$ Hence inequality in Corollary \ref{coroPoincaremodgauss}  for $f = g \circ (Id,\Psi)$ transfers to 
the same inequality for $g$ provided
$\frac{\partial f}{\partial y_i}$ is replaced  by $\partial_{i,\Psi} g.$

\bex
\mylabel{ex:unisphere}
Let $k \geq 2$ be an integer, $\SS^{k-1} \subset \RR^k$ the unit $k-1$ dimensional sphere, $\Omega = (\RR^+_{*} \times \SS^{k-1})^{\NN},$ and let  $E = (\RR^k_{*})^{\NN}.$ A typical point in $\Omega$ will be written as $(\rho,\theta) = (\rho^i,\theta^i)$ and a typical point in $E$ as $y = (y^j).$ Now consider the change of variables 
$\Psi : E \mapsto \Omega$  given by $\Psi(y) = (\Psi(y)^j)$  with 
$$\Psi^j(y) = \left(\|y^j\|^2, \frac{y^j}{\|y^j\|}\right).$$ The image of $(\gamma^k)^{\NN} = \gamma^{\NN}$ by $\Psi$ is the product measure $\tilde{\gamma}^{\NN}$
where 
$\tilde{\gamma}$ is the probability measure  on $\RR^+_{*} \times S^{k-1}$ defined by
$$ \tilde{\gamma}(dt dv) =  \frac{1}{r_k}e^{-t/2} t^{k/2 -1}\mathbf{1}_{t > 0} dt dv$$
Here $r_k = \int_0^{\infty}  e^{-t/2} t^{k/2 -1} dt$, and $dv$ stands for the
uniform probability  measure on $\SS^{k-1}$. For $g : \{0,1\}^S \times \Omega
\mapsto \RR$ with $g \circ (Id,\Psi) \in H_1^2(\lambda\otimes\gamma^{\NN})$, let 
$$\left (\frac{\partial g}{\partial \rho^j}(x,\rho, \theta), \nabla_{\theta^j}
  g(x,\rho,\theta) \right ) \in \RR \times T_{\theta^i} \SS^{k-1} \subset \RR
\times \RR^k$$ denote the partial gradient of $g$ with respect to the variable
$(\rho^i, \theta^i)$ where $T_{\theta^i} S^{k-1} \subset \RR^k$ stands for the tangent space of $\SS^{k-1}$ at $\theta_i.$
It is not hard to verify that for all $i \in \NN$ and $j \in \{1, \ldots, k\}$,
\begin{equation}
\label{eq:gradsphere}
\partial_{i,j,\Psi} g (x,\rho, \theta) 
= 2 \frac{\partial g}{\partial \rho^i}(x,\rho, \theta) \sqrt{\rho^i} \theta^i_j 
+ \frac{1}{\sqrt{\rho^i}} [\nabla_{\theta^i} g(x,\rho,\theta)]_j.
\end{equation}
As a consequence, we may recover in this way Corollary
\ref{coroPoincaremodgamma} when the parameter $a$ equals $k/2-1$, with $k$ an
integer. Indeed, this 
follows from (\ref{eq:gradsphere}) applied to the map $(x,\rho,\theta)
\rightarrow  g(\frac{\rho}{2\alpha})$. Concentrating on the angular part instead of the radial one, we obtain the
following modified Poincaré inequality on the sphere.
\eex


\begin{coro}[Uniform distribution on $\SS^{n}$]
Let $dv_n$ denote the normalized Riemannian probability measure on $\SS^{n} \subset \RR^{n+1}.$ For $g \in H^1_2(dv_n)$  and $i = 1, \ldots, n+1$ let $\nabla_i g(\theta)$ denote the $i^{th}$ component of $\nabla g(\theta)$ in $\RR^{n+1}$ (we see $T_{\theta} \SS^{n}$ as the vector space of $\RR^{n+1}$ consisting of vector that are orthogonal to $\theta$). Then, for $n \geq 2,$
$$$$
\begin{equation}
\label{eq:sphere}
\Var(g)\log\frac{\Var(g)}{\sum_{i=1}^\NN\left\|\Delta_i g\right\|_{1}^2}\leq
\frac{1}{n-1}
\sum_{i\in\NN}\EE\left(\left(\nabla_ig\right)^2\right)\;.
\end{equation}
\end{coro}
\begin{dem}
follows from (\ref{eq:gradsphere}) applied to the map $(x,\rho,\theta) \rightarrow  g(\theta).$ Details are left to the reader.
\end{dem}
\bex
\mylabel{ex:changementvar}
If one wants to get a result similar to Corollary
\ref{coroPoincaremodgauss} with $\gamma$
replaced by another probability measure $\nu$, one may of course perform the
usual change of variables through inverse of repartition function. In the
sequel, we denote by
\begin{equation}
\label{defg}
g(x)=\frac{1}{\sqrt{2\pi}}e^{-\frac{x^2}{2}}
\end{equation}
the density of the normalized Gaussian distribution, and by
\begin{equation}
\label{defG}
G(x)=\int_{-\infty }^xg(u)\;du
\end{equation}
its repartition function. For any function $\phi$ from $\RR$ to $\RR$, we
shall note $\tilde{\phi}$ the function from $\RR^\NN$
to $\RR^\NN$ such that $(\tilde{\phi}(x))_j=\phi(x_j)$.
\begin{coro}[Unidimensional change of variables]
\label{corochgtvar}
 Let $\nu$ be a probability on  $\RR^+$  absolutely continuous with respect to the  Lebesgue
measure, with density $h$ and repartition function
$$H(t)=\int_0^th(u)\;du\;.$$
Let $\{0,1\}^S\times\RR^n$ be equipped with the probability measure
$\lambda_p^S\otimes\nu^{\otimes \NN}$. Then, for every function $f$ on $\{0,1\}^S\times\RR^n$ 
such that $f\circ (Id, \widetilde{H^{-1}\circ G})\in H_1^2\left(\lambda_p^S\otimes \gamma^\NN\right)$,
$$\Var(f)\log\frac{\Var(f)}{\sum_{s\in S}\left\|\Delta_s f\right\|_{1}^2+\sum_{i\in\NN}\left\|\Delta_i f\right\|_{1}^2}\leq
c_{LS}(p)\sum_{s\in S}\EE\left((\Delta_s
  f)^2\right)+2\sum_{i\in\NN}\EE\left(\left(\nabla_i f\right)^2\right)\;,$$
where for every integer $i$,
$$\nabla_if(x,y)=\psi(y_i)\frac{\partial
  f}{\partial y_i}(x,y)\;,$$
and $\psi$ is defined on $I=\{t\geq 0\mbox{ s.t. }h(t)>0\}$:
$$\forall t\in I,\;\psi(t)=\frac{g\circ G^{-1}(H(t))}{h(t)}\;.$$
\end{coro}
\begin{dem}
It is a straightforward consequence of Corollary
\ref{coroPoincaremodgauss}, applied to $f\circ (Id, \widetilde{H^{-1}\circ G})$.\end{dem}
\eex


\section{A general exponential concentration inequality}
\mylabel{sec:concentration}

In this section, we show how one can deduce from Proposition \ref{proppoincaremod} an exponential concentration
inequality for a function $F$ of independent variables. We shall prove in section
\ref{sec:FPP}, in the context of First Passage Percolation, that this new general concentration inequality may in certain cases
improve on the ones due to Talagrand
\cite{Talagrand94b,Talagrand95,Talagrand96a}, or Boucheron et
al. \cite{Boucheronetal03}. The reason why we can get stronger results is that
Proposition \ref{proppoincaremod} is generally stronger than a simple Poincaré
inequality, and it is well known
(see Ledoux \cite{Ledoux96}, Corollary 3.2 p.49 and Theorem 3.3 p.50) that a Poincaré inequality for a measure $\mu$
implies an exponential concentration inequality for any Lipschitz function of
a random variable with distribution $\mu$. This can be achieved through
applying the Poincaré inequality to $\exp (\theta f)$, and then performing
some recurrence. This last step is essentially contained
in the following simple version, adapted to our case, of Corollary 3.2 p.49 in
\cite{Ledoux96}.
\begin{lemm}
\label{lemconcentration}
Let $f$ be a measurable real function on a probability space
$(\mathbf{X},\mathcal{X},\mu)$, and $K$ a positive constant. Suppose that for
any real number $\theta< \frac{1}{2\sqrt{K}}$, the function $x\mapsto
e^{\theta f(x)}$ is in $L^1(\mu)$, and:
$$\Var(e^{\frac{\theta f}{2}})\leq K\theta^2\EE(e^{\theta f})\;.$$
Then,
$$\forall t\geq 0,\;\mu\left(f- \int f\;d\mu>t\sqrt{K}\right)\leq 4e^{-t}\;,$$
and,
$$\forall t\geq 0,\;\mu\left(f- \int f\;d\mu<-t\sqrt{K}\right)\leq 4e^{-t}\;.$$
\end{lemm}

Now, we can state our general concentration inequality. For any function $F$ on a product space
$\left(\mathbb{X}_i,\mathcal{X}_i,\mu_i\right)_{i\in I}$, we define the
following quantities, which play an important role in Theorem
\ref{thmconcentrationgenerale} (The notation is that of section \ref{sec:logsobfkbis}).
$$W_{i,+}(x)=\int\left(F(x^{-i},y_i)-F(x)\right)_+\;d\mu_i(y_i)\;,$$
where $h_+=\sup\{h,0\}$.
$$W_+(x)=\sum_{i\in I}W_{i,+}\;.$$
Remark that similar quantities are involved in the work of Boucheron
et al. \cite{Boucheronetal03,Boucheronetal05}.

\begin{theo}
\label{thmconcentrationgenerale}
 Let $\left(\mathbb{X}_i,\mathcal{X}_i,\mu_i\right)_{i\in I}$,
 $(\mathcal{A}_i)_{i\in I}$ and $(R_i)_{i\in I}$ be as in Proposition
 \ref{proppoincaremod}, and satisfying all the hypotheses therein. Let $F$ be a function in $L^2(\Pi_{i\in I}\mathbb{X}_i)$. Define
$$r=\sup_{i\in I}\sqrt{\EE(W_{i,+}^2)},\;$$
$$s=\sqrt{\EE\left(W_+^2\right)}\;.$$
Define, for every real number $K>ers$:
$$l(K)=\frac{K}{\log\frac{K}{rs\log\frac{K}{rs}}}\;.$$
Suppose that there exists a real number $K>ers$ such that,
for every $\theta$ such that 
$|\theta|\leq \frac{1}{2\sqrt{l(K)}}$, $e^{\frac{\theta}{2} F}$ belongs to $\mathcal{A}^I$, and:
\begin{equation}
\label{eqhypmajoenergie}
\sum_{i\in I}\EE(R_i(e^{\frac{\theta}{2}F}))\leq K\theta^2\EE(e^{\theta
  F})\;.
\end{equation}
Then, denoting $\mu=\otimes_{i\in I}\mu_i$, for every $t>0$:
$$\mu(F-\EE(F)\geq t\sqrt{l(K)})\leq 4e^{-t}\;,$$
$$\mu(F-\EE(F)\leq -t\sqrt{l(K)})\leq 4e^{-t}\;.$$
\end{theo}
\begin{dem}
For any function $f$ in $L^1\left(\Pi_{i\in I}\mathbb{X}_i\right)$,  any $i\in I$, $x\in\Pi_{i\in I}\mathbb{X}_i$ and $y_i\in\mathbb{X}_i$,
\begin{eqnarray*}
\left\|\Delta_i f\right\|_1 &= &\int |\int( f(x)- f(x^{-i},y_i))
d\mu_i(y_i)|\;d\mu(x)\;,\\
&\leq &\int \int| f(x)- f(x^{-i},y_i)
|\;d\mu(x)d\mu_i(y_i)\;,\\
&= &2\int \int (f(x)- f(x^{-i},y_i))_+
\;d\mu(x)d\mu_i(y_i)\;,\\
&= &2\int \int (f(x)- f(x^{-i},y_i))_-
\;d\mu(x)d\mu_i(y_i)\;,\\
\end{eqnarray*}
where $h_-=\sup\{-h,0\}$. On the other hand,
\begin{eqnarray*}
\left(e^{\frac{\theta}{2}F(x^{-i},y_i)}-
  e^{\frac{\theta}{2}F(x)}\right)_+
&=&e^{\frac{\theta}{2}F(x)}\left(e^{\frac{\theta}{2}(F(x^{-i},y_i)-F(x))}-1\right)_+\;,\\
&\leq
  &e^{\frac{\theta}{2}F(x)}\left(\frac{\theta}{2}(F(x^{-i},y_i)(x)-F(x))\right)_+\;,\\
&=
  &\left\lbrace\begin{array}{l}\frac{|\theta|}{2}e^{\frac{\theta}{2}F(x)}\left(F(x^{-i},y_i)-F(x)\right)_+\mbox{
  if }\theta>0\;,\\ \frac{|\theta|}{2}e^{\frac{\theta}{2}F(x)}\left(F(x^{-i},y_i)-F(x)\right)_-\mbox{
  if }\theta<0\;.\end{array}\right.
\end{eqnarray*}
Therefore,
$$\left(e^{\frac{\theta}{2}F(x^{-i},y_i)}-
  e^{\frac{\theta}{2}F(x)}\right)_+\leq \left\lbrace\begin{array}{l}\frac{|\theta|}{2}e^{\frac{\theta}{2}F(x)}\left(F(x^{-i},y_i)-F(x)\right)_+\mbox{
  if }\theta>0\;,\\ \frac{|\theta|}{2}e^{\frac{\theta}{2}F(x^{-i},y_i)}\left(F(x)-F(x^{-i},y_i)\right)_+\mbox{
  if }\theta<0\;.\end{array}\right.$$
Since $F(x)$ and $F(x^{-i},y_i)$ have the same distribution under
$\mu\otimes\mu_i$, we get, for any real number $\theta$,
$$\left\|\Delta_ie^{\frac{\theta  F}{2}}\right\|_1\leq
|\theta|\int \int\left(F(x^{-i},y_i)-F(x)\right)_+\;d\mu_i(y_i)\;e^{\frac{\theta}{2}F(x)}\;d\mu(x)\;.$$
And, using Cauchy Schwarz inequality,
$$\sum_{i\in I}\left\|\Delta_ie^{\frac{\theta
        F}{2}}\right\|_1\leq
        |\theta|\sqrt{\EE\left(W_+^2\right)\EE\left(e^{\theta F}\right)}\;.$$
But we also have, again using Cauchy Schwarz inequality,
$$\left\|\Delta_ie^{\frac{\theta
      F}{2}}\right\|_1\leq|\theta|\sqrt{\EE\left(W_{i,+}^2\right)\EE\left(e^{\theta F}\right)}\;.$$
Therefore,
\begin{equation}
\label{eqmajoenergie2sum}
\sum_{i\in I}\left\|\Delta_ie^{\frac{\theta
        F}{2}}\right\|_1^2\leq \theta^2rs\EE\left(e^{\theta F}\right)\;.
\end{equation}
Inequality (\ref{eqmajoenergie2sum}), the Poincaré inequality for $e^{\theta F}$ (Proposition
\ref{proppoincaremod}) and hypothesis (\ref{eqhypmajoenergie}) imply
that:
\begin{equation}
\label{eqbeforealter}
\forall |\theta|\leq\frac{1}{2\sqrt{l(K)}}
,\;\Var(e^{\frac{\theta}{2}F})\log\frac{\Var(e^{\frac{\theta}{2}F})}{\theta^2rs\EE\left(e^{\theta
      F}\right)}\leq K\theta^2\EE\left(e^{\theta      F}\right)\;.
\end{equation}
Therefore, we are left in front of the following alternative:

$\bullet$ either $\Var(e^{\frac{\theta
        F}{2}})\leq \theta^2\frac{K}{\log\frac{K}{rs}}\EE(e^{\theta F})$,

$\bullet$ or $\Var(e^{\frac{\theta F}{2}})>
        \theta^2\frac{K}{\log\frac{K}{rs}}\EE(e^{\theta  F })$. But in
        this case, plugging this minoration into the logarithm of
        inequality (\ref{eqbeforealter}) leads to:
$$\Var(e^{\frac{\theta F}{2}})\leq  \theta^2\frac{K}{\log\frac{K}{rs\log\frac{K}{rs}}}\EE(e^{\theta \tilde{f}})\;.$$
In any case, for any $|\theta|\leq \frac{1}{2\sqrt{l(K)}}$,
$$\Var(e^{\frac{\theta
        F}{2}})\leq  \theta^2 l(K)\EE(e^{\theta \tilde{f}})\;.$$
The result follows from Lemma \ref{lemconcentration}.
\end{dem}
\begin{rem}
It is well known, through Herbst's argument (see e.g. \cite{Ledoux01}
Theorem 5.3 p.95), that condition
(\ref{eqhypmajoenergie}) implies a subexponential concentration inequality of
the form:
$$\mu(|F-\EE(F)|\geq 2t\sqrt{K})\leq 2e^{-t^2}\;.$$
In the applications to follow, $K/rs$ is big, and therefore $l(K)$ is small
compared to $K$. Therefore, at the price of trading the subgaussian
behaviour against a subexponential one, Theorem \ref{thmconcentrationgenerale} shows that when $K/rs$
is big, the fluctuations of $F$ are lower than $\sqrt{l(K)}$, which is small
compared to $\sqrt{K}$.
\end{rem}

Let us give a closer look at the case where for every $i$, $\mu_i$ is the invariant measure of a diffusion process
with carré du champ $\Gamma_i$.  Naturally associated with this diffusion process, a Sobolev logarithmic
inequality for $\mu_i$ has the form (if it exists):
$$\Ent_{\mu_i}(f^2)\leq c_i\EE_{\mu_i}(\Gamma_i(f,f))\;,$$
where $c_i$ is a positive constant. Furthermore, we have the following property:
$$\Gamma_i(\Phi(f),g)=\Phi'(f)\Gamma_i(f,g)\;,$$
which leads to:
$$R_i(e^{\frac{\theta}{2} F})^2=c_i\Gamma_i(e^{\frac{\theta}{2}
  F},e^{\frac{\theta}{2} F})=c_i\frac{\theta^2}{4}e^{\theta
  F}\Gamma_i(F,F)\;.$$
Therefore, condition (\ref{eqhypmajoenergie}) becomes:
$$\EE\left(e^{\frac{\theta}{2}
  F}\sum_{i\in I}c_i\Gamma_i(F,F)\right)\leq 4K\EE(e^{\theta
  F})\;.$$
The main work to satisfy  condition (\ref{eqhypmajoenergie}) is to
  bound from below, and somewhat independently from $e^{\frac{\theta}{2}
  F}$, the quantity $\sum_{i\in I}\Gamma_i(F,F)$. In some particular cases,
  and notably percolation, this quantity is upperbounded by $F$ itself. And it
  is possible to show, following Boucheron et al. \cite{Boucheronetal03},
 that, at least for small $\theta$, $\EE(Fe^{\theta
  F})$ is upper bounded by a constant times $\EE(F)\EE(e^{\theta
  F})$. More generally, one can state the following result.

\begin{coro}
\label{coroconcentrationgenerale}
 Let $\left(\mathbb{X}_i,\mathcal{X}_i,\mu_i\right)_{i\in I}$,
 $(\mathcal{A}_i)_{i\in I}$ and $(R_i)_{i\in I}$ be as in Proposition
 \ref{proppoincaremod}, and satisfying all the hypotheses therein. Let $F$ be a function in $L^2(\Pi_{i\in I}\mathbb{X}_i)$. Define
$$r=\sup_{i\in I}\sqrt{\EE(W_{i,+}^2)},\;$$
$$s=\sqrt{\EE\left(W_+^2\right)}\;.$$
Define, for every real number $K>ers$:
$$l(K)=\frac{K}{\log\frac{K}{rs\log\frac{K}{rs}}}\;.$$
Suppose that there exists two constants $C$ and $D$ such that, denoting:
$$A_{CD}=4C\EE(F)+D\left(1+\frac{2}{C}\right)\;,$$
we have
\bdes
\iti $C\leq \sqrt{l\left(A_{CD}\right)}$,
\itii $A_{CD}4C\EE(F)+D\left(1+\frac{2}{C}\right)\geq ers$,
\itiii for every $\theta$ such that
$|\theta|\leq \frac{1}{2\sqrt{l\left(A_{CD}\right)}}$, $e^{\theta F}$ belongs to $\mathcal{A}^I$, and:
\begin{equation}
\label{eqhypmajoenergiebis}
\sum_{i\in I}\EE(R_i(e^{\frac{\theta}{2}F}))\leq C\theta^2\EE(Fe^{\theta
  F})+D\theta^2\EE(e^{\theta
  F})\;.
\end{equation}
\edes
Then, denoting $\mu=\otimes_{i\in I}\mu_i$, for every $t>0$:
$$\mu(F-\EE(F)\geq t\sqrt{l\left(A_{CD}\right)})\leq 4e^{-t}\;,$$
$$\mu(F-\EE(F)\leq -t\sqrt{l\left(A_{CD}\right)})\leq 4e^{-t}\;,$$
\end{coro}
\begin{dem}
The only thing to prove is that condition (\ref{eqhypmajoenergie}) holds with $K=4C\EE(F)+D(1+\frac{2}{C})$. This will
follow from condition (\ref{eqhypmajoenergiebis}) and a variation on the theme of Herbst's argument due to Boucheron
et al. \cite{Boucheronetal03}. Indeed, recall that using the tensorisation of entropy, the
logarithmic Sobolev inequalities for each $\mu_i$ imply that:
$$\forall g\in \mathcal{A}^I,\; \Ent_{\mu}(g^2)\leq
\sum_{i\in I}\EE_\mu\left(R_i(g)^2\right)\;.$$
Let us apply this inequality to $g=e^{\frac{\theta}{2} F}$, and use condition
(\ref{eqhypmajoenergiebis}). For every $\theta$ such that $|\theta|\leq \frac{1}{2\sqrt{l(4C\EE(F))}}$,
$$\Ent_{\mu}\left(e^{\theta F}\right)\leq C\theta^2\EE(Fe^{\theta
  F})\;.$$
This may be written as:
\begin{equation}
\label{eqinterMassart}
\theta\EE(Fe^{\theta F})-\EE(e^{\theta F})\log\EE(e^{\theta F})\leq C\theta^2\EE(Fe^{\theta
  F})+D\theta^2\EE(e^{\theta
  F})\;.
\end{equation}
First, suppose that $\theta$ is positive. The proof of Theorem 5 in \cite{Boucheronetal03} shows that, for every
  $\theta <\frac{1}{C}$,
$$\log \EE(e^{\theta F})\leq  \frac{\theta}{1-\theta C}\EE(F)+D\frac{\theta^2}{1-\theta C}\;,$$
and equation (\ref{eqinterMassart}) implies that, for every
  $\theta <\frac{1}{C}$,
$$\EE(Fe^{\theta F})\leq \frac{\EE(F)+D\theta}{(1-\theta C)^2}\EE(e^{\theta
  F})\;.$$
If $\theta$ is negative, $e^{\theta
  F}$ is decreasing in $F$, and it follows from Chebyshev's association inequality
that (see e.g. \cite{HardyLittlewoodPolya52} p.43):
$$\EE(Fe^{\theta F})\leq \EE(F)\EE(e^{\theta F})\;.$$
Now, we gather the case where $\theta$ is positive and the case where it is
negative. Condition {\bf (i)} implies that $ \frac{1}{2C} \geq \frac{1}{2\sqrt{l(4C\EE(F))}}$,
and therefore, for every $\theta$ such that $|\theta|\leq \frac{1}{2\sqrt{l(4C\EE(F))}}$,
\begin{eqnarray*}
\sum_{i\in I}\EE(R_i(e^{\frac{\theta}{2}F}))&\leq& C\theta^2\EE(Fe^{\theta
  F})+D\theta^2\EE(e^{\theta
  F})\\
&\leq& \left(4C\EE(F)+D\left(1+\frac{2}{C}\right)\right)\theta^2\EE(e^{\theta
  F})\;,
\end{eqnarray*}
and the result follows from Theorem \ref{thmconcentrationgenerale}.
\end{dem}
The
main lesson that we can remember from Corollary \ref{thmconcentrationgenerale} is the following (very) informal
statement. 
\vspace*{3mm}

\begin{center}{\it 
If $F$ is a lipschitz function of a large number of variables, each
of which contributes at most to an amount $\delta$, then $F$ has fluctuations
of order $O(\sqrt{\EE(F)/log \frac{1}{\delta}})$, and there is an exponential
control for these fluctuations.}
\end{center}
\vspace*{5mm}

\section{Application to First Passage Percolation}
\mylabel{sec:FPP}

\subsection{Continuous edge-times distributions}
\mylabel{subsec:continuousFPP}
It turns out that Corollary \ref{coroconcentrationgenerale} is particularly well suited
to adapt the argument of Benjamini, Kalai and Schramm
\cite{BenjaminiKalaiSchramm03} to show that the passage time from the origin
to a vertex $v$ satisfies
an exponential  concentration inequality at the rate $O(\sqrt{|v|/\log |v|})$
when the edges have a $\Gamma(a,b)$ distribution with $a\geq 1/2$. This includes
the important case of exponential distribution, for which First Passage
Percolation becomes equivalent to a version of Eden growth model (see for instance Kesten
\cite{Kesten84}, p.130). We do not want
to restrict ourselves to those distributions. Nevertheless, due to the particular strategy that we
adopt, we can only prove our result for some continuous edge times distributions
which behave roughly like a gamma distribution. Please note that the
definition given below differs (one assumption is removed) from the definition
of a nearly gamma distribution that was stated in the preliminary paper \cite{BenaimRossignolarxiv06}.

\begin{defi}
\mylabel{defneargamma}
{\rm
 Let $\nu$ be a probability on  $\RR^+$  absolutely continuous with respect to the  Lebesgue
measure, with density $h$ and repartition function
$$H(t)=\int_0^th(u)\;du\;.$$
Define:
$$I=\{t \geq 0\mbox{ such that }h(t) > 0\}\;,$$
and $\psi : I \mapsto \RR$ the map:
$$\psi(y)=\frac{g\circ G^{-1}(H(y))}{h(y)}\;.$$
Let $A$ be a positive real number. The probability measure $\nu$ will be said
to be \emph{nearly gamma} provided it satisfies the following set of conditions:
\bdes
\iti $I$ is an interval;
\itii $h$ restricted to $I$ is continuous;
\itiii There exists a positive real number $A$ such that
$$\forall y\in I,\;\psi(y)\leq
A\sqrt{y}\;.$$
\edes
}
\end{defi}
If we want to emphasize the dependance on $A$ in the above definition, we
shall say that $\nu$ is \emph{nearly gamma with bound $A$}. In Definition \ref{defneargamma}, condition $(iii)$ is of course the most
tedious to check.  A simple sufficient condition for a
probability measure to be nearly gamma will be given in Lemma
\ref{lemneargamma}, the proof of which relies on the following asymptotics for
the Gaussian repartition function $G$.
\begin{lemm}
\label{lemasymptotgG}As $x$ tends to $-\infty$,
$$G(x)=g(x)\left(\frac{1}{|x|}+o\left(\frac{1}{x}\right)\right)\;,$$
and as $x$ tends to $+\infty$,
$$G(x)=1-g(x)\left(\frac{1}{x}+o\left(\frac{1}{x}\right)\right)\;.$$
Consequently,
$$g\circ G^{-1}(y)\stackrel{y\rightarrow 0}{\sim}y\sqrt{-2\log y}\;,$$
and
$$g\circ G^{-1}(y)\stackrel{y\rightarrow 1}{\sim}(1-y)\sqrt{-2\log (1-y)}\;.$$
\end{lemm}
\begin{dem}
A simple change of variable $u=x-t$ in $G$ gives:
$$G(x)=g(x)\int_{0}^{+\infty}e^{-\frac{t^2}{2}+xt}\;dt\;,$$
Integrating by parts, we get:
\begin{eqnarray*}
G(x)&=&g(x)\left(-\frac{1}{x}+\frac{1}{x}\int_{0}^{+\infty}te^{-\frac{t^2}{2}+xt}\;dt\right)\;,\\
&=&g(x)\left(\frac{1}{|x|}+o\left(\frac{1}{x}\right)\right)\;,\\
\end{eqnarray*}
as $x$ goes to $-\infty.$ Since
$G(-x)=1-G(x)$, we get that, as $x$ goes to $+\infty$:
$$G(x)=1-g(x)\left(\frac{1}{x}+o\left(\frac{1}{x}\right)\right)\;.$$
Let us turn to the asymptotic of $g\circ G^{-1}(y)$ as $y$ tends to zero. Let
$x=G^{-1}(y)$, so that ``$y$ tends to zero'' is equivalent to ``$x$ tends to
$-\infty$''. One has therefore,
\begin{eqnarray*}
G(x)&=&\frac{g(x)}{|x|}(1+o(1))\;,\\
\log G(x)&=&\log g(x)-\log |x| +o(1)\;,\\
&=&-\frac{x^2}{2}-\log |x| +O(1)\;,\\
\log G(x)&=&-\frac{x^2}{2}(1+o(1))\;,\\
|x|&=&\sqrt{-2\log G(x)}\;.\\
\end{eqnarray*}
Since $g(x)=|x|G(x)(1+o(1))$,
$$g(x)=G(x)\sqrt{-2\log G(x)}(1+o(1))\;,$$
and therefore,
$$g\circ G^{-1}(y)=y\sqrt{-2\log y}(1+o(1))\;,$$
as $y$ tends to zero. The asymptotic of $g\circ
G^{-1}(y)$ as $y$ tends to 1 is derived in the same way.
\end{dem}
Given two functions $r$ and $l$, we write $l(x) = \Theta(r(x))$ as $x$ goes to $x*$ provided there exist positive constants $C_1 \leq C_2$ such that
$$C_1 \leq \liminf_{x \rightarrow x*}  \frac{r(x)}{l(x)} \leq \limsup_{x \rightarrow x*} \frac{r(x)}{l(x)} \leq C_2.$$

\begin{lemm}
\label{lemneargamma}
Assume that condition $(i)$ and $(ii)$ of  Definition \ref{defneargamma} hold. Let $0 \leq \underline{\nu} < \overline{\nu} \leq \infty$ denote the endpoints of $I.$ Assume furthermore condition $(iii)$ is replaced by conditions $(iv)$ and $(v)$ below.
\bdes
\itiv There exists $\alpha >-1$ such that as $x$ goes to $\underline{\nu}$,
$$h(x)=\Theta\left((x-\underline{\nu})^\alpha\right)\;,$$
\itv $\overline{\nu} < \infty$ and there exists $\beta >-1$ such that as $x$ goes to $\overline{\nu}$,
$$h(x)=\Theta\left((\overline{\nu}-x)^\beta\right)\;,$$
\noindent or $\overline{\nu} = \infty$ and
$$\exists A>\underline{\nu},\;\forall t\geq A,\; C_1h(t)\leq \int_t^{\infty} h(u)\;du\leq C_2h(t)\;,$$
\noindent where $C_1$ and $C_2$ are positive constants.
\edes
Then, $\nu$ is nearly gamma.
\end{lemm}
\begin{dem}
Since $h$ is a continuous function on $]\underline{\nu},\overline{\nu}[$, it
attains its minimum on every compact set included in
$]\underline{\nu},\overline{\nu}[$. The minimum of $h$ on
$[a,b]$ is therefore strictly positive as
soon as $\underline{\nu}<a\leq b<\overline{\nu}$. In order to show
that condition $(iii)$ holds, we thus have to concentrate on the behaviour of the
function $\psi$ near $\underline{\nu}$ and $\overline{\nu}$. Condition $(iv)$
implies that, as $x$ goes to $\underline{\nu}$,
\begin{equation}
\label{DLHgauche}
H(x)=\Theta\left((x-\underline{\nu})^{\alpha+1}\right)\;.
\end{equation}
This, via Lemma \ref{lemasymptotgG}, leads to
\begin{equation}
\label{DLpsigauche}
\psi(x)=\Theta\left((x-\underline{\nu})\sqrt{-\log(x-\underline{\nu})}
  \right)\;,
\end{equation}
as $x$ goes to $\underline{\nu}$. Similarly, if $\overline{\nu} < \infty$,
  condition $(v)$ implies that, as $x$ goes to $\overline{\nu}$,
\begin{equation}
\label{DLHdroite}
H(x)=\Theta\left((\overline{\nu}-x)^{\beta+1}\right)\;,
\end{equation}
which leads via Lemma \ref{lemasymptotgG} to
\begin{equation}
\label{DLpsidroite}
\psi(x)=\Theta\left((\overline{\nu}-x)\sqrt{-\log(\overline{\nu}-x)}
  \right)\;,
\end{equation}
as $x$ goes to $\overline{\nu}$. Therefore, if $\overline{\nu} < \infty$,
   condition $(iii)$ holds.

Now, suppose that $\overline{\nu}=\infty$. Condition $(v)$ implies:
$$\forall t\geq A,\;1/C_2\leq \frac{h(t)}{\int_t^{\infty} h(u)}\;du\leq 1/C_1\;.$$
Integrating this inequality between A and $y$ leads to the existence of three
positive constants $B$, $C'_1$ and $C'_2$ such that:
$$\forall y\geq B,\;C'_1y\leq \log\frac{1}{1-H(y)}\leq C'2y\;.$$
Thus,
\begin{equation}
\label{DLpsiinfini}
\forall y\geq B,\;C_1\sqrt{C'_1y}\leq \psi(y)\leq C_2\sqrt{C'_2y}\;.
\end{equation}
This, combined with equation (\ref{DLpsigauche}) proves that condition $(iii)$
holds and concludes the proof of Lemma \ref{lemneargamma}.
\end{dem}
\brem
With the help of Lemma \ref{lemneargamma}, it is easy to check that most usual
distributions are nearly gamma. This includes all gamma and beta
distributions, as well as any probability measure whose
density is bounded away from 0 on its support, and notably the uniform
distribution on
$[a,b]$, with $0\leq a<b$. Nevertheless, remark that some distributions which 
have a sub-exponential upper tail may not satisfy the assumptions of Lemma
\ref{lemneargamma}, and be nearly gamma, though. For example, this is the case of the
distribution of $|N|$, where $N$ is a standard Gaussian random variable.
\erem
Now, we can state the main result of this article.
\begin{theo}
\label{theoFPPexpo} Let $\nu$ be a nearly gamma probability measure with an
exponential moment, i.e we suppose that there exists $\delta>0$ such that:
$$\int e^{\delta x}\;d\nu(x)<\infty\;.$$
Let $\mu$ denote the measure $\nu^{\otimes E}$. Then, there exist two positive
constants $C_1$ and $C_2$ such that, for any $|v|\geq 2$, and any positive real number $t\leq |v|$,
$$\mu\left(|d_x(0,v)-\int
d_x(0,v)\;d\mu(x)|> t\sqrt{\frac{|v|}{\log |v|}}\right)\leq C_1e^{-C_2t}\;.$$
\end{theo}
\begin{dem}
What we present here borrows many ideas from Kesten \cite{Kesten93} and of course Benjamini et
al. \cite{BenjaminiKalaiSchramm03}. We would like to apply Corollary \ref{coroconcentrationgenerale}  to the
function $e^{\theta f_v}$, for $\theta\leq \sqrt{\frac{\log |v|}{|v|}}$. In fact, we will be able to use Corollary \ref{coroconcentrationgenerale}, but not exactly for $f_v$,
  and not exactly for any nearly gamma distribution. The first step is indeed
  to work with a version of $\nu$ with bounded support. Precisely, we shall use the following lemma which is an easy
adaptation of Kesten's Lemma 1, p.309 in \cite{Kesten93}. 
\begin{lemm}
\label{lemtroncature}
Let $\nu$ be a nearly gamma distribution with bound $A$. Suppose that $\nu$ admits an exponential moment, i.e there exists $\delta>0$ such that:
$$\int e^{\delta x}\;d\nu(x)<\infty\;.$$
Then there exists a sequence of probability measures $(\nu_k)_{k\geq 2}$, positive constants $C_3,C_4,C_5$
and a positive integer $k_{\nu}$ with the following properties:
\bdes
\iti For every $k$, the support of $\nu_k$ is included in $[0,C_5\log k]$,
\itii If $k\geq k_{\nu}$, $\nu_k$ is a nearly gamma distribution with bound $A$.
\itiii If $k=|v|$ and $k\geq 2$, for every $t$ greater than $2C_3\sqrt{\frac{\log
    |v|}{|v|}}$, 
\begin{eqnarray*}
&&\nu\left(|d_x(0,v)-\EE_\nu(d_x(0,v))|>t\sqrt{\frac{|v|}{\log
      |v|}}\right)\\
 &\leq &3e^{-C_3|v|}+C_4e^{-\frac{\gamma}{8}t\sqrt{\frac{|v|}{\log
      |v|}}} +\tilde{\nu}\left(|d_x(0,v)-\EE_{\tilde{\nu}}(d_x(0,v))|>\frac{t}{4}\sqrt{\frac{|v|}{\log
      |v|}}\right) \;.
\end{eqnarray*}
\itiv If $k\geq 2$, $\nu_k$ is stochastically smaller than $\nu_{k+1}$ and $\nu$.
\edes
\end{lemm}
\begin{dem}
Kesten's argument in \cite{Kesten93} is simply to consider the truncated
edge times at $C_5\log |v|$. We cannot use this directly  because we have to deal with continuous
distribution. Instead, we can repatriate the mass beyond $2C_5\log |v|$, and
spread it continuously over $[C_5\log |v|,2C_5\log |v|]$. This mass is small, of
course. Precisely, thanks to the exponential moment assumption, for every
positive number $c$,
$$\nu([c\log |v|,+\infty[)\leq \int e^{\delta x}\;d\nu(x)\frac{1}{|v|^{\delta c}}\;.$$
Let $u$ be a continuous density on the real line with support included in
$[0,1]$ and $C_5$ a positive constant to be fixed later. We define
  $\nu_k$ to be the continuous distribution on the real line with
  density:
$$\forall x\in\RR,h_k(x)=\left(h(x)+u\left(\frac{x-C_5\log k}{C_5\log
      k}\right)\frac{\nu([2C_5\log k,+\infty[)}{C_5\log k}\right)\II_{x\leq
      2C_5\log k}\;.$$
Statements {\bf(i)} and {\bf(iv)} are obvious. To see that {\bf(ii)} holds,  let $H_k$ be the
repartition function of $\nu_k$. Obviously,
\begin{eqnarray*}
\forall x\leq C_5\log|v|,\;h_k(x)&=& h(x)\;,\\
\forall x\leq 2C_5\log|v|,\;h_k(x)&\geq& h(x)\;,\\
\forall x\geq C_5c\log|v|,\;h_k(x)&=&0\;,
\end{eqnarray*}
and therefore,
\begin{eqnarray*}
\forall x\leq C_5\log|v|,\;H_k(x)&=& H(x)\;,\\
\forall x\in \RR,\;H_k(x)&\geq& H(x)\;,\\
\forall x\geq 2C_5\log|v|,\;H_k(x)&=&1\;.
\end{eqnarray*}
Observe now that $g\circ G^{-1}$ is decreasing on $[1/2,1]$, and that
$$\forall x\geq C_5\log k,\;H_k(x)\geq H(x)\geq 1-\int e^{\delta x}\;d\nu(x)\frac{1}{k^{\delta C_5}}\;.$$
Therefore, let $k_{\nu}=\left\lceil\left(2\int e^{\delta
    x}\;d\nu(x)\right)^{\frac{1}{\delta C_5}}\right\rceil$,
$$\forall k\geq k_{\nu},\;\forall x\leq
    2C_5\log k,\;\frac{g\circ G^{-1}(H_k(x))}{h_k(x)}\leq
    \psi(x)\;.$$
This implies that the distributions $(\nu_k)_{k\geq k_\nu}$ are all nearly
    gamma with  the same bound $A$.

It remains to prove  {\bf(iv)}. We define the following coupling $\pi_k$ of $(\nu,\nu_k)$:
$$\int g(x,y)\;d\pi_k(x,y)=\int g(x,H_k^{-1}(H(x)))\;d\nu(x)\;.$$
Denote by $\gamma=\gamma(\tilde{x})$ the ($\nu_k^{\otimes E}$ a.s unique) $\tilde{x}$-geodesic from 0 to
$v$. The following inequalities hold for $\pi_k$-almost every $(x,\tilde{x})$.
\begin{eqnarray*}
0\leq d_x(0,v)-d_{\widetilde{x}}(0,v)&\leq&
\sum_{e\in\gamma}x_e-\sum_{e\in\gamma}\widetilde{x}_e\;,\\
&\leq&
\sum_{e\in\gamma}x_e\II_{x_e>C_5\log k}\;.\\
\end{eqnarray*}
Now, if $k=|v|$, we choose to take $C_5=\frac{4d}{\delta}$, and the end of the
proof follows exactly Kesten's Lemma 1 p.~309 in
\cite{Kesten93}.
\end{dem}

Now, we suppose that $k=|v|\geq k_\nu$ and we shall work with $\nu_k$, whose
support is included in $[0,2C_5\log |v|]$. Let us define $\mu_k=\nu_k^{\otimes
  E}$. In the whole proof, $Y$ shall denote
a random variable with distribution $\nu$. Remark that, thanks to
part {\bf(iv)} of Lemma \ref{lemtroncature}
$$\int e^{\delta x}\;d\nu_k(x)\leq \int e^{\delta x}\;d\nu(x)=\EE(e^{\delta Y})\;,$$
and, for any positive real number $\alpha$,
$$\int x^\alpha\;d\nu_k(x)\leq \int  x^\alpha\;d\nu(x)=\EE(Y^\alpha)\;.$$
A crucial idea in the work of Benjamini Kalai and Schramm is to work with a randomised version of $f_v$ in order to take
a full benefit of Corollary \ref{corochgtvar}. This randomisation trick relies on the following lemma:
\begin{lemm}
\label{lemaveraging}
There exists a constant $c>0$, such that, for every $m\in \NN^*$, there exists a
function $g_m$ from $\{0,1\}^{m^2}$ to $\{0,\ldots,m\}$ such that:
$$\max_{y\in\{0,\ldots,m\}}\lambda(x\mbox{ s.t. }g_m(x)=y)\leq\frac{c}{m}\;,$$
and
$$\forall q\in\{1,\ldots ,m^2\},\nabla_qg_m\in\{0,1\}\;,$$
where 
$$\nabla_qg(x)=g(x_1,\ldots,x_{q-1},1,x_{q+1},\ldots,x_{m^2})-g(x_1,\ldots,x_{q-1},0,x_{q+1},\ldots,x_{m^2})\;.$$
\end{lemm}
Since Benjamini et al. do not give a full proof for this lemma, we
offer the following one.
\begin{dem}
From Stirling's Formula,
$$\binom{m^2}{\lfloor m^2/2\rfloor}.
\frac{m}{2^{m^2}}\xrightarrow[n\rightarrow \infty]{}\frac{1}{\sqrt{2\pi}}\;,$$
and this implies that the following supremum is finite:
$$c_1=\sup\left\lbrace 2\binom{m^2}{\lfloor
      m^2/2\rfloor}.\frac{m}{2^{m^2}}\mbox{ s.t. }m\in\NN^*\right\rbrace\;.$$
Notice also that $c_1\geq 1$. Now, let $\preceq$ denote the alphabetical order
      $\{0,1\}^{m^2}$, and let us list the elements in $\{0,1\}^{m^2}$ as
      follows:
$$(0,0,\ldots ,0)=x_1\preceq x_1\preceq \ldots \preceq x_{2^{m^2}}=(1,1,\ldots ,1)\;.$$
For any $m$ in $\NN^*$, we define the
      following integer:
$$k(m)=\left\lceil\frac{2^{m^2}}{m}\right\rceil\;,$$
and the following function on $\{0,1\}^{m^2}$:
$$\forall i\in \{1,\ldots
,2^{m^2}\},\;g_m(x_i)=\left\lfloor\frac{i}{k(m)}\right\rfloor\;.$$
Remark that $g_m(x_{2^{m^2}})\leq m/c_1\leq 1$. Therefore, $g$ is a function
from $\{0,1\}^{m^2}$ to $\{0,\ldots m\}$. Now, suppose that $x_i$ and $x_l$
differ from exactly one coordinate. Then, 
\begin{eqnarray*}
|i-l|&\leq & \binom{m^2}{i}+\binom{m^2}{l}\;,\\
&\leq & 2\binom{m^2}{\lfloor m^2/2\rfloor}\;,\\
&\leq & c_1\frac{2^{m^2}}{m}\;,\\
&\leq & k(m)\;.
\end{eqnarray*}
Consequently,
\begin{eqnarray*}
g_m(x_i)-g_m(x_l)&\leq&\left\lfloor\frac{l}{k(m)}+1\right\rfloor-\left\lfloor\frac{l}{k(m)}\right\rfloor\;,\\
&=& 1\;,
\end{eqnarray*}
which implies that $\nabla_qg_m\in\{0,1\}$. Finally, for any $y\in\{0,\ldots
m\}$, $g$ takes the value $y$ at most $k(m)$ times, and
\begin{eqnarray*}
\lambda(x\mbox{ s.t. }g_m(x)=y)& \leq & \frac{k(m)}{2^{m^2}}\;,\\
&\leq &\frac{c_1}{m}+\frac{1}{2^{m^2}}\;,\\
&\leq &\frac{2c_1}{m}\;.
\end{eqnarray*}
So the lemma holds with $c=2c_1$.
\end{dem}
Now, we define our randomised version of $f_v$ as follows. Let $m$ be a
positive integer, to be fixed later, and $S=\{1,\ldots,d\}\times\{1,\ldots,m^2\}$. Let
$c>0$ and $g_m$ be as in Lemma \ref{lemaveraging}. As in 
\cite{BenjaminiKalaiSchramm03}, for any $a=(a_{i,j})_{(i,j)\in S}\in \{0,1\}^S$, let
$$z=z(a)=\sum_{i=1}^dg_m(a_{i,1},\ldots ,a_{i,m^2})\mathbf{e_i}\;,$$
where $(\mathbf{e_1},\ldots ,\mathbf{e_d})$ denotes the standard basis of
$\ZZ^d$. We now equip the space $\{0,1\}^S\times \RR_+^E$ with the probability
measure $\lambda\otimes \mu_k$, where $\lambda:=\lambda_{\frac{1}{2}}^S$ is the uniform
measure on $\{0,1\}^S$, and we define the following function $\tilde{f}$ on $\{0,1\}^S\times
\RR_+^E$:
$$\forall (a,x)\in \{0,1\}^S\times \RR_+^E, \;
\tilde{f}(a,x)=d_x(z(a),v+z(a))\;.$$
When $m$ is not too big, $f$ and $\tilde{f}$ are not too
far apart. 
\begin{lemm}
\label{lemfftilde}
For any positive real number $t$,
$$\mu_k(|f-\EE(f)|>t)\leq  \lambda\otimes\mu_k(|\tilde{f}-\EE(\tilde{f})|>t/2)+e^{-\frac{\delta t}{4m}}\EE\left(e^{\delta Y}\right)\;.$$
\end{lemm}
\begin{dem}
Let $\alpha(a)$ be a path from 0 to $z(a)$,
such that $|\alpha(a)|=|z(a)|$ (here, $|\alpha|$ is the number of edges in
$\alpha$). Let $\beta(a)$ denote a path disjoint from $\alpha(a)$, which goes
from $v$ to $v+z(a)$. Then,
\begin{eqnarray*}
|\tilde{f}(a,x)-f(x)|&\leq& d_x(0,z(a))+d_x(v,v+z(a))\;,\\
&\leq &\sum_{e\in\alpha (a)}x_e+\sum_{e\in\beta (a)}x_e\;,\\
\end{eqnarray*}
which is stochastically dominated by a sum of $2m$ independent variables
$Y_1,\ldots ,Y_{2m}$ with distribution $\nu$. Remark that, due to the translation invariance of the distribution of $f$ under $\mu_k$, $f$ and $\tilde{f}$ have the same mean against $\lambda\otimes\mu_k$. Thus, using $|z|\leq m$, we have:
$$\mu_k(|f-\EE(f)|>t)\leq \lambda\otimes \mu_k(|\tilde{f}-\EE(\tilde{f})|>t/2)+\lambda\otimes \mu_k (|f-\tilde{f}|>t/2)\;.$$
Now, by Markov's inequality, we get that for any positive real number $t$, 
\begin{eqnarray*}
\lambda\otimes
\mu_k (|f-\tilde{f}|>t/2)&\leq&\PP\left(\sum_{i=1}^{2m}Y_i>\frac{t}{2}\right)\;,\\
&=&\PP\left(\frac{\delta}{2m}\sum_{i=1}^{2m}Y_i>\frac{\delta t}{4m}\right)\;,\\
&\leq&e^{-\frac{\delta t}{4m}}\EE\left(e^{\frac{\delta}{2m} Y}\right)^{2m}\;,\\
&\leq&e^{-\frac{\delta t}{4m}}\EE\left(e^{\delta Y}\right)\;.\\
\end{eqnarray*}
This concludes the proof of this lemma.
\end{dem}

It remains to bound $\lambda\otimes\mu_k(|\tilde{f}-\EE(\tilde{f})|>t)$. To this end, we will
use an adaptation of Corollary \ref{coroconcentrationgenerale}, applied to
$F=\tilde{f}$. Denote, for any $s$ in $S$ and any $e$ in $E$,
$$W_{s,+}=\int\left(F(x^{-s},y_s)-F(x)\right)_+\;d\beta_{1/2}(y_s)\;,$$
and
$$W_{S,+}=\sum_{s\in S}W_{s,+}\;,$$
$$W_{e,+}=\int\left(F(x^{-e},y_e)-F(x)\right)_+\;d\nu_k(y_e)\;,$$
and
$$W_{E,+}=\sum_{e\in E}W_{e,+}\;.$$
Applying Corollary \ref{corochgtvar} with $p=1/2$ (note that $c_{LS}(1/2)=2$),
we can get the following minor adaptation of Corollary
\ref{coroconcentrationgenerale}. The notations are those of Corollary \ref{corochgtvar} and Definition
\ref{defneargamma}.

\begin{prop}
\label{propadaptFPP}
 Let $\nu$ be a probability on  $\RR^+$  absolutely continuous with respect to the  Lebesgue
measure, with density $h$ and repartition function
$$H(t)=\int_0^th(u)\;du\;.$$
Let $\{0,1\}^S\times\RR^n$ be equipped with the probability measure
$\lambda_p^S\otimes\nu^{\otimes \NN}$. Let $F$ be a function from $\{0,1\}^S\times\RR^n$ to $\RR$.  Define
$$r_S=\sup_{s\in S}\sqrt{\EE(W_{s,+}^2)}\;,$$
$$s_S=\sqrt{\EE\left(W_{S,+}^2\right)}\;,$$
$$r_E=\sup_{e\in E}\sqrt{\EE(W_{e,+}^2)}\;,$$
$$s_E=\sqrt{\EE\left(W_{E,+}^2\right)}\;,$$
and
$$K_{ES}=r_Ss_S+r_Es_E\;.$$
Define, for every real number $K>eK_{ES}$:
$$l(K)=\frac{K}{\log\frac{K}{K_{ES}\log\frac{K}{K_{ES}}}}\;.$$
Suppose that there exists three positive real numbers $C,$ $D$ and $A_{CD}$ such that:
\bdes
\iti $C\leq \sqrt{l(A_{CD})}$,
\itii $A_{CD}\geq \sup\{eK_{ES},4C\EE(F)+D\left(1+\frac{2}{C}\right)\}$,
\itiii for every $\theta$ such that
$|\theta|\leq \frac{1}{2\sqrt{l(A_{CD})}}$,  
 $e^{\theta F}\circ (Id, \widetilde{H^{-1}\circ G})\in
 H_1^2\left(\lambda_p^S\otimes \gamma^\NN\right)$ and:
\begin{equation}
\label{eqhypmajoenergieS}
\sum_{s\in S}\left\|\Delta_s(e^{\frac{\theta}{2}F})\right\|_2^2\leq D\theta^2\EE(e^{\theta
  F})\;,
\end{equation}
and:
\begin{equation}
\label{eqhypmajoenergieE}
\sum_{e\in E}\left\|\nabla_e(e^{\frac{\theta}{2}F})\right\|_2^2\leq C\theta^2\EE(Fe^{\theta
  F})\;,
\end{equation}
where for every  $e$ in $E$,
$$\nabla_ef(x,y)=\psi(y_e)\frac{\partial
  f}{\partial y_e}(x,y)\;,$$
and $\psi$ is defined on $I=\{t\geq 0\mbox{ s.t. }h(t)>0\}$:
$$\forall t\in I,\;\psi(t)=\frac{g\circ G^{-1}(H(t))}{h(t)}\;.$$
\edes
Then, denoting $\mu=\lambda_S\otimes\gamma^E$, for every $t>0$:
$$\mu(F-\EE(F)\geq t\sqrt{l(A_{CD})})\leq 4e^{-t}\;,$$
$$\mu(F-\EE(F)\leq -t\sqrt{l(A_{CD})})\leq 4e^{-t}\;.$$
\end{prop}
First, we need to prove that $e^{\theta\tilde{f}}\circ
(Id,\widetilde{H^{-1}\circ G})$ belongs to
$H_1^2(\lambda\otimes\gamma^\NN)$ when $\nu$ is nearly gamma. This is the aim of the following Lemma.
\begin{lemm}
\label{lemderiveefv}
 If $\nu$ is nearly gamma, and has bounded support, for any positive number $\theta$, the function $e^{\theta f_v}\circ\widetilde{H^{-1}\circ G}$ belongs to
  $H_1^2(\gamma^\NN)$, $e^{\theta \tilde{f}}\circ
(Id,\widetilde{H^{-1}\circ G})$ belongs to
  $H_1^2(\lambda\otimes\gamma^\NN)$. Furthermore, conditionally to $z$, there is almost surely only one $x$-geodesic from $z$ to
$z+v$, denoted by $\gamma_x(z)$, and:
$$\frac{\partial
  \tilde{f}}{\partial x_e}(a,x)=\II_{e\in\gamma_x(z(a))}\;.$$
\end{lemm}
\begin{dem}
The fact that $e^{\theta f_v}\circ\widetilde{H^{-1}\circ G}$ and $e^{\theta \tilde{f}}\circ
(Id,\widetilde{H^{-1}\circ G})$ are in $L^2$ is obvious since $\nu$ has bounded support. We shall prove that $e^{\theta f_v}\circ\widetilde{H^{-1}\circ G}$ satisfies the integration by part formula {\bf
  (a)} of the definition of $H_1^2$. The similar result for $e^{\theta \tilde{f}}$ is
obtained in the same way. Now, we fix $x^{-e}$ in $(\RR^+)^{E(\ZZ^d)\setminus\{e\}}$. We denote by $g_e$ the function
defined on $\RR^+$ by:
$$g(y)=f_v(x^{-e},y)\;.$$
We will show that there is a nonnegative real number $y_\infty$ such that:
\begin{equation}
\label{claimg}
\left\lbrace\begin{array}{l}\forall y\leq y_\infty,\;g(y)=g(0)+y\\
\mbox{ and }\\
\forall y > y_\infty,\;g(y)=g(y_\infty)\end{array}\right.
\end{equation}
For any $n\geq |v|$, let us denote by $\Gamma_n$ the set of paths from 0 to $v$ whose number of
edges is not greater than $n$. We have:
$$g(y)=\inf_{n\geq |v|}g_n(y)\;,$$
where 
$$g_n(y)=\inf_{\gamma\in\Gamma_n}\sum_{e'\in \gamma}(x^{-e},y)_{e'}\;.$$
The functions $g_n$ form a nonincreasing sequence of nondecreasing functions:
$$\forall n\geq|v|,\;\forall y\in\RR^+,\;\forall y'\geq y,\;g_{n+1}(y)\leq
g_n(y)\leq g_n(y')\;.$$
In particular, this implies that for every $y$ in $\RR^+$,
$$g(y)=\lim_{n\rightarrow \infty}g_n(y)\;.$$ 
Now, we claim that, for every $n\geq |v|+3$, there exists $y_n\in\RR^+$ such that:
\begin{equation}
\label{claimgn1}
\left\lbrace\begin{array}{l}\forall y\leq y_n,\;g_n(y)=g_n(0)+y\\
\mbox{ and } \\
\forall y > y_n,\;g_n(y)=g_n(y_n)\end{array}\right.
\end{equation}
and furthermore, 
\begin{equation}
\label{claimgn2}
\mbox{ the sequence }(y_n)_{n\geq|v|+3}\mbox{ is nonincreasing.}
\end{equation}
Indeed, since
$\Gamma_n$ is a finite set, the infimum in the definition of $g_n$
is attained. Let us call a path which attains this infimum an
$(n,y)$-geodesic and let $\tilde{\Gamma}(n,y,e)$ be the set of $(n,y)$-geodesics
which contain the edge $e$. Remark that as soon
as $n\geq |v|+3$, there
exists a real number $A$ such that $e$ does not belong to any $(n,A)$-geodesic:
it is enough to take $A$ greater than the sum of the length of three edges forming a path between the
end-points of the edge $e$. Therefore, the following supremum is finite:
$$y_n=\sup\{y\in\RR^+\mbox{ s.t. }\tilde{\Gamma}(n,y,e)\not =\emptyset\}\;.$$
Now, if $e$ belongs to an $(n,y)$-geodesic $\gamma$, for any
$y'\leq y$, $\gamma$ is an $(n,y')$-geodesic to which $e$ belongs, and $g_n(y)-g_n(y')=y-y'$. If
$\tilde{\Gamma}(n,y,e)$ is empty, then for any $y'\geq y$, $e$ does
not belong to any $(n,y')$-geodesic, and $g_n(y)=g_n(y')$. This proves that:
$$\forall y< y_n,\;g_n(y)=g_n(0)+y\;,$$
$$\forall y,y' > y_n,\;g_n(y)=g_n(y')\;.$$
Since $g_n$ is continuous, we have proved claim (\ref{claimgn1}). Now remark
that if $e$ does not belong to any $(n,y)$-geodesic, then $e$ does not belong to
any $(n+1,y)$-geodesic, since $\Gamma_n\subset\Gamma_{n+1}$. Therefore,
$y_{n+1}\leq y_n$, and this proves claim (\ref{claimgn2}). Since
$(y_n)_{n\geq|v|+3}$ is nonnegative, it converges to a nonnegative number
$y_\infty$ as $n$ tends to infinity. Now, let $n$ be a integer greater than
$|v|+3$:
$$\forall n\geq N,\;\forall y,y'>y_n,\;g_n(y)=g_n(y')\;.$$
Since $y_n\leq y_N$, 
$$\forall n\geq N,\;\forall y,y'>y_N,\;g_n(y)=g_n(y')\;.$$
Letting $n$ tend to infinity in the last equation, we get:
$$\forall N\geq |v|+3,\;\forall y,y'>y_N,\;g(y)=g(y')\;.$$
Therefore,
$$\forall y,y'>y_\infty,\;g(y)=g(y')\;.$$
On the other side,
$$\forall n\geq |v|+3,\;\forall y\leq y_n,\;g_n(y)=g_n(0)+y\;.$$
Since $y_n\geq y_\infty$,
$$\forall n\geq |v|+3,\;\forall y\leq y_\infty,\;g_n(y)=g_n(0)+y\;.$$
Letting $n$ tend to infinity in the last expression, we get:
$$\forall y\leq y_\infty,\;g(y)=g(0)+y\;.$$
Finally, $g$ is continuous. Indeed, the convergent sequence $(g_n)$ is
uniformly equicontinuous, since all these functions are 1-Lipschitz, and the continuity of
$g$ follows from Arzelà-Ascoli Theorem. We have proved claim
(\ref{claimg}). Remark that $y_\infty=y_\infty(x^{-e})$ depends on
$x^{-e}$. We define, for any $x^{-e}$,
$$h_e(x^{-e},x_e)=\left\lbrace\begin{array}{l}1\mbox{ if } x_e\leq y_\infty(x^{-e}) \\ 0\mbox{ if } x_e >y_\infty(x^{-e})\end{array}\right.\;.$$
It is easy to see that, for any
smooth function $F:\RR\rightarrow \RR$ having compact support, for any
$x^{-e}$,
\begin{equation}
\label{eqIPP}
-\int_{\RR} F'(x_e)e^{\theta f_v(x^{-e},x_e)}\;dx_e=\theta\int_{\RR}
F(x_e)h_e(x^{-e},x_e)e^{\theta f_v(x^{-e},x_e)}\;dx_e\;.
\end{equation}
It is known that there is almost surely a geodesic from 0 to $v$ (see
\cite{Howard04} for instance), i.e the
infimum in the definition of $f_v$ is attained with probability
1. Furthermore, in this setting, where the distribution of the lengths is
continuous, there is almost surely only one unique $x$-geodesic from $0$ to
$v$. For any  $x$, we shall denote by $\gamma_x(0)$ the unique $x$-geodesics from
  $0$ to $0+v$. Then, with $\nu$-probability 1, one can see from the definitions
  of $y_n$ and $y_\infty$ that:
\begin{equation}
\label{eqexpressionder}h_e(x^{-e},x_e)=\II_{e\in\gamma_x(0)}\;.
\end{equation}
Performing the change of variable $x\mapsto \widetilde{H^{-1}\circ G}$ in
equation (\ref{eqIPP}), one gets the integration by parts formula {\bf (a)}
for $e^{\theta f_v}\circ \widetilde{H^{-1}\circ G}$, with the following partial derivative with
respect to $x_e$:
$$x\mapsto \theta\psi(x_e)h_e(\widetilde{H^{-1}\circ G}(x))e^{\theta f_v}\;.$$
The expression of $\frac{\partial
  \tilde{f}}{\partial x_e}(a,x)$ is derived in the same way than
(\ref{eqexpressionder}).
\end{dem}

Now, we want to apply Proposition \ref{propadaptFPP} to $F=\tilde{f}$.

{\bf Bound on $\sum_{s\in
    S}\left\|\Delta_s(e^{\frac{\theta}{2}F})\right\|_2^2$}

Here, we can perform a quite rough upper bound,
        since there are not many elements in $S$. For any $a\in \{0,1\}^S$, and any $q$ in $S$, denote by $\tau_qa$ the element
of $\{0,1\}^S$ obtained from $a$ by flipping the coordinate $q$. Then, for any
function $g$ on $\{0,1\}^S$,
\begin{eqnarray*}
\left\|\Delta_q e^{\theta g/2}
    \right\|_{p}^p&=& \frac{1}{4}\sum_{a\in\{0,1\}^S}\left|e^{\frac{\theta}{2}
    g(a)}-e^{\frac{\theta}{2} g(\tau_qa)}\right|^p\lambda(a)\;,\\
&=&\frac{1}{2}\sum_{a:\;\theta g(a)>\theta g(\tau_qa)}e^{\frac{\theta p}{2}
    g(a)}\left(1-e^{\frac{\theta}{2} (g(\tau_qa)-g(a))}\right)^p\lambda(a)\;,\\
&\leq&\frac{|\theta|^p}{2^{p+1}}\sum_{a:\;\theta g(a)>\theta
    g(\tau_qa)}e^{\frac{\theta p}{2}
    g(a)}|g(a)-g(\tau_qa)|^p\lambda(a)\;,\\
&\leq&\frac{|\theta|^p}{2^{p+1}}\sum_{a\in\{0,1\}^S}e^{\frac{\theta p}{2}
    g(a)}|g(a)-g(\tau_qa)|^p\lambda(a)\;.\\
&=&\frac{|\theta|^p}{2}\left\|e^{\theta g/2}\Delta_q g
    \right\|_{p}^p\;.\\
\end{eqnarray*}
 According to Lemma \ref{lemaveraging}, for any $q\in \{0,1\}^{m^2}$,
$\nabla_qg_m\in\{0,1\}$. Therefore, for any $s=(i,q)\in S$, 
$$|\Delta_s \tilde{f}|\leq\frac{1}{2}(
x_{(z,z+\mathbf{e_1})}+x_{(z+v,z+v+\mathbf{e_1})})\;.$$
Therefore we get the following bounds:
\begin{equation}
\label{majoenergieS}
\sum_{s\in S}\left\|\Delta_se^{\frac{\theta
        \tilde{f}}{2}}\right\|_2^2\leq \theta^2C_5^2m(\log|v|)^2\EE(e^{\theta\tilde{f}})\;.
\end{equation}

{\bf Bound on $r_S$}

\begin{equation}
\label{majosupnorm1S}
r_S\leq 2\sqrt{\EE(Y^2)}
\end{equation}

{\bf Bound on $s_S$}

\begin{equation}
\label{majosommenorm1S}
r_S\leq 2m
\end{equation}

{\bf Bound on $\sum_{e\in E}\left\|\nabla_e(e^{\frac{\theta}{2}F})\right\|_2^2$}

Let $A$ be as
  in Definition \ref{defneargamma}. Since $\nu_k$ is nearly gamma with bound
  $A$ (see Lemma \ref{lemtroncature}),
$$\sum_{e\in E}\left\|e^{\theta \tilde{f}/2}\nabla_e\tilde{f}\right\|_2^2\leq A\EE\left(\tilde{f}e^{\theta \tilde{f}}\right)\;.$$

{\bf Bound on $s_E$}
Remark that:
$$(\tilde{f}(x^{-e},y_e)-\tilde{f}(x))_+\leq y_e\II_{e\in\gamma_x(z)}\;,$$
and $\gamma_x(z)$ is independent from $y_e$. Therefore,
\begin{equation}
\label{eqmajograddiscretE}
0\leq W_{e,+}\leq \EE(Y)\II_{e\in\gamma_x(z)}\;,
\end{equation}
which leads to:
$$0\leq W_{E,+}\leq \EE(Y)|\gamma_x(z)|\;,$$
and:
$$s_E\leq \EE(Y)\sqrt{\EE(|\gamma_x(z)|^2)}\;.$$

Now, following Kesten \cite{Kesten93}, p.308, we claim that there exists some
        constant $C_6$, depending only on $\nu$ (and not on $k$) such that:
\begin{equation}
\label{eqnorm1geodesiccarre}
\EE_{\nu_k}\left(|\gamma_x(z)|^2\right)\leq C_6|v|^2\;.
\end{equation}
Indeed, for any $a>0$ and $y>0$,
\begin{eqnarray*}
\mu_k(|\gamma_x(0)|\geq y|v|)&\leq &\mu_k(f_v\geq ay|v|)\\
&&+\mu_k(\exists \mbox{ a
  self-avoiding path }r\mbox{ starting at }0\mbox{ of at least }\\
&&y|v|\mbox{
  steps but with }\sum_{e\in r}x_e<ay|v|)\;.
\end{eqnarray*}
Proposition 5.8 of Kesten \cite{Kesten84} shows that for a suitable
  $a>0$, the second term in the right-hand side of the above inequality is at
  most $Ce^{-C'y|v|}$ for some constants $C$ and $C'$. Further more,
  $a$, $C$ and $C'$ do not depend on $k$: it suffices to choose them for
  $\nu_{k_\nu}$, and the same constants work for any $k\geq k_\nu$ (see part
  {\bf(iv)} of Lemma
  \ref{lemtroncature} and the remark of Kesten \cite{Kesten93} p.309). On the
  other hand, $f_v$ is dominated by the sum of $|v|$ independent variables
  with distribution $\nu$, $X_1,\ldots ,X_{|v|}$. Thus,
\begin{eqnarray*}
\EE\left(|\gamma_x(z)|^2\right)&=&\EE\left(|\gamma_x(z)|^2\right)\;,\\
&=&|v|^2\int_0^\infty\mu(|\gamma_x(0)|^2>y|v|^2)\;dy\;,\\
&\leq &|v|^2\int_0^\infty\mu\left(\left(\sum_{i=1}^{|v|}X_i\right)^2\geq
  ay|v|^2\right)+|v|^2C\int_0^\infty e^{-C'\sqrt{y}|v|}\;dy\;,\\
&=&\frac{1}{a^2}\EE\left(\left(\sum_{i=1}^{|v|}X_i\right)^2\right)+2C\int_0^\infty
  te^{-C't}\;dt\;,\\
&\leq&C_6|v|^2\;.\\
\end{eqnarray*}
This proves claim (\ref{eqnorm1geodesiccarre}).Therefore, 
$$s_E\leq \sqrt{C_6}\EE(Y)|v|\;.$$

{\bf Bound on $r_E$}

From inequality (\ref{eqmajograddiscretE}), we get:
$$r_E\leq\EE(Y)\sqrt{\sup_{e\in E}\PP(e\in\gamma_x(z))}\;.$$
Now, we use the fact that for any fixed $z$, $\mu$ is invariant under translation by $z$.
\begin{eqnarray*}
\PP(e\in\gamma_x(z))&=&
\EE_{\lambda}\left(\EE_{\mu}\left(\II_{e-z\in\gamma_{x}(0)}\right)\right)\;,\\
&=&\EE_{\mu}\left(\sum_{e'\in\gamma_{x}(0)}\EE_{\lambda}\left(\II_{e-z=e'}\right)\right)\;,\\
&=&\EE_{\mu}\left(\sum_{e'\in\gamma_{x}(0)}\PP_{\lambda}(z=e-e')\right)\;,\\
&\leq& \sup_{z_0}\PP(z=z_0)\EE_{\mu}\left(|\gamma_{x}(0)\cap
  \mathcal{Q}_e|\right)\;,\\
&\leq& \sup_{z_0}\PP(z=z_0)\EE_{\mu}\left(|\gamma_{x}(0)\cap \mathcal{B}_e|\right)\;,\\
\end{eqnarray*}
where $\mathcal{Q}_e=\left\lbrace e'\in E\left(\ZZ^d\right)\mbox{
    s.t. }\PP(z=e-e')>0\right\rbrace \subset \mathcal{B}_e=e+\mathcal{B}(0,dm)$. Using Lemma \ref{lemaveraging}, 
$$ \sup_{z_0}\PP(z=z_0)\leq \left(\frac{c}{m}\right)^d\;.$$
Now, we claim that
\begin{equation}
\label{eqsupnorm1geodesiccarre}
\EE_{\mu_k}\left(|\gamma_{x}(0)\cap \mathcal{B}_e|\right)\leq C_7m^{d-1}\;,
\end{equation}
We proceed as we did to obtain (\ref{eqnorm1geodesiccarre}). Indeed, for any $a>0$ and $y>0$,
\begin{eqnarray*}
\mu_k\left(|\gamma_{x}(0)\cap \mathcal{B}_e|\geq ym\right)&\leq &\mu_k\left(\sum_{e'\in\gamma_{x}(0)\cap \mathcal{B}_e}x_{e'}\geq ay|v|\right)\\
&&+\sum_{w\in\partial\mathcal{B}_e}\mu_k(\exists \mbox{ a
  self-avoiding path }r\mbox{ starting at }\\
&&w\mbox{ of at least }ym\mbox{  steps but with }\sum_{e\in r}x_e<aym)\;.
\end{eqnarray*}
We use again the constants $a$, $C$ and $C'$ arising from Proposition 5.8 of Kesten
  \cite{Kesten84}, and which depend on $\nu$, but not on $k$. Remark that
  there are at most $(dm)^{d-1}$ vertices in $\partial\mathcal{B}_e$. On the
  other hand, let $r$ be a deterministic path going through
every vertex of the surface of the ball $\mathcal{B}_e$, and such that
there is a constant $C''$ (depending only on $d$) such that $|r|\leq C''m^{d-1}$. From the
definition of a geodesic, we get:
$$f_v\leq\sum_{e'\in r}x_e\;.$$
Thus,
\begin{eqnarray*}
\EE\left(|\gamma_{x}(0)\cap \mathcal{B}_e|\right)&=&m\int_0^\infty\mu_k(|\gamma_x(0)\cap \mathcal{B}_e|>ym)\;dy\;,\\
&\leq &m\int_0^\infty\mu_k\left(\sum_{e'\in r}x_e\geq
  aym\right)+m(dm)^{d-1}\int_0^\infty e^{-C'ym}\;dy\;,\\
&=&\frac{1}{a}\EE_{\mu_k}\left(\sum_{e'\in r}x_e\right)+\frac{2C}{C'}(dm)^{d-1}\int_0^\infty te^{-C't}\;dt\;,\\
&\leq&C_7m^{d-1}\;.\\
\end{eqnarray*}
This proves claim  (\ref{eqsupnorm1geodesiccarre}). Therefore:
\begin{eqnarray}
\nonumber r_E&\leq&\EE(Y)\sqrt{\left(\frac{c}{m}\right)^dC_7m^{d-1}}\;,\\
\label{majorE} r_E&\leq&\frac{C_8}{m^{\frac{1}{2}}}\;.
\end{eqnarray}

{\bf End of the proof}

Now, we choose  $m=\lceil
|v|^{1/4}\rceil$. Define $C=A$, $D=C_5^2m(\log|v|)^2$. The bounds obtained before lead to:
$$K_{ES}=O(|v|^{\frac{7}{8}})\;,$$
and:
$$4C\EE(F)+D\left(1+\frac{2}{C}\right)=O(|v|)\;.$$
So we can choose $A_{CD}=C_4|v|$, with $C_4$ a positive constant, such that
{\bf (ii)}  of Proposition
\ref{propadaptFPP} applied to $F=\tilde{f}$ is satisfied. It is clear that, for $|v|$
large enough, conditions {\bf (i)} and {\bf (iii)}  are also satisfied. Remark
also that:
$$l(A_{CD})=O\left(\frac{|v|}{\log |v|}\right)\;.$$
Therefore, there exists a
constant $C_{12}$ such that for every $t>0$:
\begin{equation}
\label{eqconcentrationftildedroite}
\mu_k\left(\tilde{f}-\EE(\tilde{f})>t\sqrt{\frac{|v|}{\log |v|}}\right)\leq
        4e^{-C_{12}t}\;.
\end{equation}
and:
\begin{equation}
\label{eqconcentrationftildegauche}
\mu_k\left(\tilde{f}-\EE(\tilde{f})<-t\sqrt{\frac{|v|}{\log |v|}}\right)\leq
        4e^{-C_{12}t}\;.
\end{equation}
Lemmas \ref{lemfftilde} and \ref{lemtroncature} conclude the proof of Theorem \ref{theoFPPexpo}.
\end{dem}

\begin{rem}
Inequalities
(\ref{eqconcentrationftildedroite}), (\ref{eqconcentrationftildegauche}) and Lemma
\ref{lemfftilde} imply, after integration, that the
variance of $f_v$ is of order $O(|v|/\log |v|)$. Of course, we do not need the
assumption that $\nu$ has a bounded support to obtain such a result. Instead, we just need $\nu$ to have a second moment. The proof mimics \cite{BenjaminiKalaiSchramm03}, and the ideas presented here. Details may be found in \cite{BenaimRossignolarxiv06}, which is a preliminary  version of the present paper.
\end{rem}

\subsection{Bernoulli distributions}

The method developped in subsection \ref{subsec:continuousFPP} applies also to
the case where the edge-times are distributed according to a Bernoulli law $\nu =(1-p)\delta_a+p\delta_b$, and $a$ is strictly positive. The proof
follows exactly the same pattern as the proof of the nearly gamma case, except
that:
\begin{enumerate}
\item one does not need Lemma \ref{lemtroncature}, since $\nu$ has bounded
  support,
\item the geodesic is not almost surely unique anymore,
\item the energy $\sum_{e\in
    E}\EE\left(R_e\left(e^{\frac{\theta}{2}\tilde{f}}\right)^2\right)$ is
    different.
\end{enumerate}
Point 1 is just  good news. Point 2 is not a problem: the bounds on $s_E$,
$s_S$, $r_E$ and $r_S$ remain valid if we choose for $\gamma_x(z)$ one
geodesic among all the possible ones. So we shall only show how to circumvent
point 3, i.e how one can bound $\sum_{e\in
    E}\EE\left(R_e\left(e^{\frac{\theta}{2}\tilde{f}}\right)^2\right)$, where 
$$R_e(f)=\sqrt{c_{LS}(p)}\Delta_ef\;.$$
First, imitating the proof of Theorem \ref{thmconcentrationgenerale}, we write:
$$\sum_{e\in
    E}\EE\left(R_e\left(e^{\frac{\theta}{2}\tilde{f}}\right)^2\right)\leq
    c_{LS}(p)\frac{\theta^2}{4}\EE\left(V_{E,+}e^{\theta\tilde{f}}\right)\;,$$
where:
$$V_{E,+}=\sum_{e\in E}\int
(\tilde{f}x^{-e},y_e)-\tilde{f}(x))_+^2\;d\nu(y_e)\;.$$
Now,
\begin{eqnarray*}
V_{E,+}&\leq &\sum_{e\in E} \int (b-a)^2\II_{e\in\gamma_x(z)}\;d\nu(y_e)\;,\\
&=&(b-a)^2|\gamma_x(z)|_1\;,\\
&\leq &\frac{(b-a)^2}{a}\tilde{f}\;,
\end{eqnarray*}
Therefore,
\begin{equation}
\label{eqboundbernoullinon0}\sum_{e\in
    E}\EE\left(R_e\left(e^{\frac{\theta}{2}\tilde{f}}\right)^2\right)\leq
    c_{LS}(p)\frac{\theta^2}{4}\frac{(b-a)^2}{a}\EE\left(\tilde{f}e^{\theta\tilde{f}}\right)\;.
\end{equation}
The bound (\ref{eqboundbernoullinon0}) allows us to obtain the following equivalent
    of Theorem \ref{theoFPPexpo} in the case of Bernoulli distributions.
\begin{prop}
\label{propFPPexpoBernoulli} Let $a$ and $b$ be two real numbers such that
    $0<a<b$. We define $\nu=(1-p)\delta_a+p\delta_p$  and
    $\mu=\nu^{\otimes E}$. Then, there exist two positive
constants $C_1$ and $C_2$ such that, for any $|v|\geq 2$, and any positive real number $t$,
$$\mu\left(|d_x(0,v)-\int
d_x(0,v)\;d\mu(x)|> t\sqrt{\frac{|v|}{\log |v|}}\right)\leq C_1e^{-C_2t}\;.$$
\end{prop}

\begin{rem}
{When $a =0$, the previous argument does not work, and it is hard to compare
$V_{E,+}$ to $\tilde{f}$ itself. Although the quantity $V_{E,+}$ may be
controlled when $1-p<p_c(\ZZ^d)$ via Kesten's work (see
Proposition 5.8 in  \cite{Kesten84}), we do not know how to adapt the entire 
proof to this case. 
}

\end{rem}

\section*{Acknowlegdgements}

R. Rossignol would like to thank warmly Stéphane Boucheron and Pascal Massart for
having insisted on the possible interest of \cite{Rossignol06} in the context of
First Passage Percolation.



\end{document}